\documentclass[11pt,reqno]{amsart}
\usepackage{mathrsfs}
\usepackage{amsmath}
\usepackage{amssymb, amsmath}
\usepackage{color}

\numberwithin{equation}{section}

\let\al=\alpha

\let\e=\varepsilon

\let\r=\rho

\let\f=\frac

\let\G= \Gamma

\let\Om=\Omega

\let\na=\nabla

\let\pa=\partial

\def\cS{{\mathcal S}}

\def\R{\mathbf R}
\def\Z{\mathbf Z}

\def\no{\noindent}

\def\dv{\mbox{div}}

\def\curl{\mathop{\rm curl}\nolimits}
\def\Ds{\langle D\rangle^s}
\def\eqdef{\buildrel\hbox{\footnotesize def}\over =}
\def\ef{\hphantom{MM}\hfill\llap{$\square$}\goodbreak}

\newcommand{\beq}{\begin{equation}}
\newcommand{\eeq}{\end{equation}}
\newcommand{\ben}{\begin{eqnarray}}
\newcommand{\een}{\end{eqnarray}}
\newcommand{\beno}{\begin{eqnarray*}}
\newcommand{\eeno}{\end{eqnarray*}}

\newtheorem{theorem}{Theorem}[section]
\newtheorem{definition}[theorem]{Definition}
\newtheorem{lemma}[theorem]{Lemma}
\newtheorem{proposition}[theorem]{Proposition}

\newtheorem{remark}[theorem]{Remark}
\newtheorem{Theorem}{Theorem}[section]

\newtheorem{Lemma}[Theorem]{Lemma}

\begin{document}

\title[Break-down criterion for the water-wave equation]
{Break-down criterion for the water-wave equation}

\author{Chao Wang}

\address{Institut de Math\'{e}matiques de Jussieu, CNRS UMR 7586, Universit\'{e} Paris 7,  75013 Paris and Beijing International Center For Mathematical Research,  Peking University, 100871, P. R. China}
\email{wangc@math.jussieu.fr}

\author{Zhifei Zhang}
\address{School of Mathematical Sciences, Peking University, 100871, P. R. China}
\email{zfzhang@math.pku.edu.cn}


\date{March 25, 2013}



\begin{abstract}
We study the break-down mechanism of smooth solution for the gravity water-wave equation of infinite depth. It is proved that if the mean curvature $\kappa$ of the free surface $\Sigma_t$, the trace $(V,B)$ of the velocity at the free surface,
and the outer normal derivative $\frac {\pa P} {\pa \textbf{n}}$ of the pressure $P$ satisfy
\beno
&&\displaystyle\sup_{t\in [0,T]}\|\kappa(t)\|_{L^p\cap L^2}+\int_0^T\|(\na V, \na B)(t)\|_{L^\infty}^6dt<+\infty,\\
&&\displaystyle\inf_{(t,x,y)\in [0,T]\times \Sigma_t}-\frac {\pa P} {\pa \textbf{n}}(t,x,y)\ge c_0,
\eeno
for some $p>2d$ and $c_0>0$, then the solution can be extended after $t=T$.
\end{abstract}

\maketitle

\section{Introduction}

\subsection{Presentation of the problem}

In this paper, we are concerned with the motion of an ideal,
incompressible, irrotational gravity fluid in a domain with free boundary of infinite depth:
\beno
\big\{(t,x,y)\in [0,T]\times\mathbf{R}^d\times\mathbf{R}: (x,y)\in \Omega_t\big\},
\eeno
where $\Omega_t$ is the fluid domain at time $t$ located by the free surface
\beno
\Sigma_t=\big\{(x,y)\in\mathbf{R}^d\times\mathbf{R}:y=\eta(t,x)\big\}.
\eeno
where $t,x,y$ denote the time variable, the horizontal and vertical spacial variables respectively.
Throughout this paper, we will use the notations:
\beno
\nabla=(\pa_{x_i})_{1\leq i\leq d},\quad \nabla_{x,y}=(\nabla,\pa_y),\quad \Delta=\sum_{1\leq i\leq d}\pa_{x_i}^2,\quad \Delta_{x,y}=\Delta+\pa_y^2.
\eeno

The motion of the fluid is described by the incompressible Euler equation
\beq\label{eq:euler}
\pa_t v+v\cdot \na_{x,y}v=-ge_{d+1}-\na_{x,y} P\qquad \hbox{in}\quad \Omega_t,\,\, t\geq 0,
\eeq
where $-ge_{d+1}=(0,\cdots,0,-g)$
denotes the acceleration of gravity and $v=(v_1,...,v_d,v_{d+1})$
denotes the velocity field. The incompressibility of the fluid is expressed by
\beq \label{eq:euler-d}
\textrm{div}\, v=0 \qquad \hbox{in}\quad \Omega_t,\,\, t\geq 0,
\eeq
and the irrotationality means that
\beq \label{eq:euler-c}
\textrm{curl}\,v=0 \qquad \hbox{in}\quad \Omega_t,\,\, t\geq 0.
\eeq

At the free surface, the boundary conditions are given by
\beq\label{eq:euler-b}
\pa_t \eta-\sqrt{1+|\na\eta|^2}v_n|_{y=\eta(t,x)}=0,\quad P(t,x,y)|_{y=\eta(t,x)}=0,
\quad \hbox{for}\quad t\ge 0, x\in\R^d,
\eeq
where $v_n=\textbf{n}_+\cdot v|_{y=\eta(t,x)}$, with $\textbf{n}_+:=\frac{1}{\sqrt{1+|\na
\eta|^2}}(-\na\eta, 1)^T$ denoting the outward normal vector to the free surface $\Sigma_t$.
The first equation of (\ref{eq:euler-b}) means that the free surface moves with the fluid.
In general, the pressure at the free surface is given by
\beno
P|_{y=\eta(t,x)}=-\kappa\na\cdot\big(\frac{\na\eta}{\sqrt{1+|\na\eta|^2}}\big)\qquad \hbox{for}\quad t\geq 0,
x\in \R^d,
\eeno
where $\kappa\ge 0$ is the surface tension coefficient. In this paper, we will consider
the case without surface tension. In such case, the pressure at the free surface can be set to zero.

As in \cite{Craig, Lan}, we use an alternative formulation of the water wave system (\ref{eq:euler})-(\ref{eq:euler-b}).
From (\ref{eq:euler-d}) and (\ref{eq:euler-c}), there exists a
potential flow function $\phi$ such that $v=\na_{x,y}\phi$ and
\beq\label{eq:vp}
\Delta_{x,y}\phi=0 \qquad \hbox{in}\quad \Omega_t,\,\, t\geq 0.
\eeq
The boundary condition (\ref{eq:euler-b}) can be expressed in terms of $\phi$
\beq \label{eq:vp-b}
\pa_t\eta-\sqrt{1+|\na\eta
|^2}\pa_{n_+} \phi|_{y=\eta(t,x)}=0, \qquad \hbox{for}\quad t>0,
\,\,x\in\R^d,
\eeq
where we denote $\pa_{n_+}\eqdef\textbf{n}_+\cdot\na_{x,y}$. The Euler's
equation (\ref{eq:euler}) can be put into Bernoulli's form
\beq\label{eq:euler-ber}
\pa_t \phi+\frac{1}{2}|\na_{x,y} \phi|^2+gy=-P \qquad \hbox{in}\quad \Omega_t,\,\, t\geq 0.
\eeq

We next reduce the system (\ref{eq:vp})-(\ref{eq:euler-ber}) to a system where all the functions are evaluated at the free surface only.
For this purpose, we introduce the trace of the velocity potential $\phi$
at the free surface
\[ \psi(t,x)\eqdef\phi(t,x,\eta(t,x)), \]
and the (rescaled) Dirichlet-Neumann operator $G(\eta)$
\beno
G(\eta)\psi\eqdef\sqrt{1+|\na\eta|^2}\pa_{n_+}\phi|_{y=\eta(t,x)}.
\eeno
Taking the trace of (\ref{eq:euler-ber}) on the free surface, the
system (\ref{eq:vp})-(\ref{eq:euler-ber}) is equivalent to the system
\beq \label{eq:euler-zak}
\left\{
\begin{array}{ll}
\pa_t \eta-G(\eta)\psi=0, \\
\pa_t \psi+g\eta+\frac{1}{2}|\na
\psi|^2-\frac{(G(\eta)\psi+\na\eta\cdot\na
\psi)^2}{2(1+|\na\eta|^2)}=0,
\end{array}\right.
\eeq
which is an evolution equation for the height of the free
surface $\eta(t,x)$ and the trace of the velocity potential on the
free surface $\psi(t,x)$.

\subsection{Main result}

Let us first recall some known results on the well-posedness of the water-wave problem.
Nalimov \cite{Nal}, Yosihara \cite{Yos} and Craig \cite{Craig} proved the local well-posedness of the 2-D water-wave equation
in the case when the motion of free surface is a
small perturbation of still water. In general,
the local well-posedness of the water wave equation without surface tension
was solved by Wu \cite{Wu1,Wu2} in the case of infinite depth. See also Ambrose and Masmoudi \cite{AM1,AM2}, where they studied the well-posedness
of the water-wave equation with surface tension and zero surface tension limit. Based on the formulation (\ref{eq:euler-zak}),
Lannes \cite{Lan} proved the local well-posedness of the water-wave equation
without surface tension in the case of finite depth; while Ming and Zhang \cite{MZ} dealt with the case with surface tension.
Recently, Alazard, Burq and Zuily \cite{ABZ1, ABZ2} proved the local well-posedness of the water-wave equation with surface tension for the low regularity initial data by
using the paradifferential operator tools and Strichartz type estimates.
We should mention some recent results \cite{CL, Lin,CS1,SZ, ZZ} concerning the local well-posedness of the rotational water-wave equation.

For small initial data, Wu \cite{Wu3} proved the almost global well-posedness of 2-D water-wave equation, and Wu \cite{Wu4} and Germain, Masmoudi and Shatah \cite{GMS}
proved the global well-posedness of 3-D water-wave equation. On the other hand, Castro, Cordoba, Ferferman, Gancedo and Lopez-Fernandez \cite{CCF1} showed that
there exists smooth initial data for the water-waves equation such that the solution overturns in finite time. See \cite{CCF2, CS2} for the splash singularity.
Wu \cite{Wu5} also construct a class of self-similar solution for the 2-D water-wave equation without the gravity.

In this paper, we will study the possible break-down mechanism of the local solution of the system (\ref{eq:euler-zak}).
For the incompressible Euler equation in the whole space, Beale, Kato and Majda \cite{BKM} showed that as long as
\beno
\int_0^T\|\curl v(t)\|_{L^\infty}dt<+\infty,
\eeno
then the solution $v$ can be extended after $t=T$. For the water-wave equation, Craig and Wayne in a survey paper \cite{CW}
propose the similar problem ``{\bf How do solutions break down?}" and state:

{\it There are several versions of this question, including `` What is the lowest exponent of a Sobolev space $H^s$ in which
one can produce an existence theorem local in time?" Or one could ask ``For which $\al$ is it true that, if one knows a priori that
$\sup_{[-T,T]}\|(\eta,\psi)\|_{C^\al}<+\infty$ and that $(\eta_0,\psi_0)\in C^\infty$, then the solution is fact $C^\infty$ over the time interval $[-T,T]$?"
$\cdots\cdots$It would be more satisfying to say that the solution fails to exist because the curvature of the surface has diverged at some point, or a
related geometrical and(or) physical statement.}

For the first version of Craig-Wayne's problem, Alazard, Burq and Zuily make the important progress in a recent work \cite{ABZ3}. To state their result,
we denote by $(V,B)$ the horizontal and vertical traces of the velocity on $\Sigma_t$, which is defined by
\ben
V\eqdef \nabla \phi|_{y=\eta},\quad B\eqdef \pa_y \phi|_{y=\eta}.
\een

\begin{theorem}\label{thm:local}\cite{ABZ3}
Let $d\geq 1, s>1+\f d2$.
Assume that the initial data $(\eta_0, \psi_0)$ satisfy
\begin{eqnarray*}
\eta_0\in H^{s+1/2}(\mathbf{R}^d),\quad \psi_0\in H^{s+\f12}(\mathbf{R}^d),\quad V_0\in H^{s}(\mathbf{R}^d), \quad B_0\in H^{s}(\mathbf{R}^d).
\end{eqnarray*}
Then there exists $T>0$ such that the system (\ref{eq:euler-zak}) with the initial data
$(\eta_0, \psi_0)$ has a unique solution $(\eta,\psi)$ satisfying
\begin{gather*}
\eta\in C([0,T];H^{s+\f12}(\mathbf{R}^d)),\quad \psi\in C([0,T];H^{s+\f12}(\mathbf{R}^d)),\\
V \in C([0,T];H^{s}(\mathbf{R}^d)), \quad B\in C([0,T];H^{s}(\mathbf{R}^d)).
\end{gather*}
\end{theorem}

\begin{remark}
In fact, the authors in \cite{ABZ3} consider the case of finite depth with rough bottom, in which case one need to impose an extra condition
on the initial data:
\beno
a(0,x)\ge c>0\quad\textrm{ for }x\in \R^d,\quad \textrm{where}\quad a(t,x)=-(\pa_y P)(t,x,\eta(t,x)),
\eeno
which is the so-called Taylor sign condition. In the case of infinite depth or finite depth with the flat bottom, this condition
automatically holds, see \cite{Lan, Wu1, Wu2} for example.
\end{remark}

The goal of this paper is to answer the second version of Craig-Wayne's problem. Our main result
is stated as follows.

\begin{theorem}\label{thm:blow-up}
Let ($\eta, \psi$) be the solution of the system (\ref{eq:euler-zak}) on $[0,T]$ stated in Theorem \ref{thm:local}. If the solution
$(\eta,\psi)$ satisfies
\beno
&&M(T)\eqdef \displaystyle\sup_{t\in [0,T]}\|\kappa(t)\|_{L^p\cap L^2}+\int_0^T\|(\na V, \na B)(t)\|_{L^\infty}^6dt<+\infty,\\
&&\textrm{TS}\eqdef\displaystyle\inf_{(t,x,y)\in [0,T]\times \Sigma_t}-\frac {\pa P} {\pa \textbf{n}}(t,x,y)\ge c_0,
\eeno
for some  $p>2d$ and $c_0>0$, then we have
\beno
\sup_{t\in [0,T]}E_s(t)\le C\big(E_s(0), M(T),T, \textrm{TS}^{-1}\big).
\eeno
Especially, the solution $(\eta,\psi)$ can be extended after $t=T$. Here $\kappa(t,x)$ is the
mean curvature of the free surface $\Sigma_t$ defined by
\beno
\kappa(t,x)\eqdef\na\cdot\Big(\frac{\na\eta(t,x)}{\sqrt{1+|\na\eta(t,x)|^2}}\Big),
\eeno
and $E_s(t)\eqdef\|(\eta,\psi)(t)\|_{H^{s+\f12}}+\|(V,B)(t)\|_{H^s}$, and $C(\cdot,\cdot,\cdot)$ is an increasing function.
\end{theorem}

\begin{remark}
Theorem \ref{thm:blow-up} implies that there are three possible blow-up mechanisms for the solution of the gravity water-wave equation:
rolling-over of the surface, blow-up of the curvature or the formation of shock.
\end{remark}

\begin{remark}
The same result should be true for the case of finite depth. In a future work, we will extend a similar result to
the rotational water-wave equation.
\end{remark}

\section{Paradifferential calculus}

In this section, we recall some results about the paradifferential calculus from \cite{Me}, see also \cite{ABZ1,ABZ3}.

\subsection{Paradifferential operators}

Let us first introduce the definition of the symbol with limited spatial smoothness. We denote $W^{k,\infty}(\R^d)$ the usual Sobolev spaces for $k\in \mathbf{N}$, and
the H\"{o}lder space with exponent $k$ for $k\in (0,1)$.
\begin{definition}
Given $\rho\in [0,1]$ and $m\in \mathbf{R}$, we denote by $\Gamma^m_\rho(\mathbf{R}^d)$ the space of locally bounded functions $a(x,\xi)$ on
$\mathbf{R}^d\times \mathbf{R}^d \backslash\{0\}$, which are $C^\infty$ with respect to $\xi$ for $\xi\neq 0$ and such that, for all $\alpha\in \mathbf{N}^d$ and all
$\xi\neq 0$, the function $x\rightarrow \pa_\xi^\al a(x.\xi)$ belongs to $W^{\rho,\infty}$ and there exists a constant $C_\al$ such that
\beno
\|\pa_\xi^\al a(\cdot,\xi)\|_{W^{\rho,\infty}}\leq C_\al(1+|\xi|)^{m-|\al|}\quad \textrm{for any} \quad |\xi|\geq \f12.
\eeno
The semi-norm of the symbol is defined by
\beno
M^m_\rho(a)\eqdef\sup_{|\alpha|\leq 3d/2+1+\rho}\sup_{|\xi|\geq1/2} \|(1+|\xi|)^{|\alpha|-m}\pa^\alpha_\xi a(\cdot,\xi)\|_{W^{\rho,\infty}}.
\eeno
\end{definition}

Given a symbol $a$, the paradifferential operator $T_a$ is defined by
\ben\label{paradiff}
\widehat{T_au}(\xi)\eqdef (2\pi)^{-d}\int \chi(\xi-\eta,\eta)\widehat{a}(\xi-\eta,\eta)\psi(\eta)\widehat{u}(\eta) d \eta,
\een
where $\widehat a(\theta,\xi)$ is the Fourier transform of $a$ with respect to the first variable;
the $\chi(\theta,\xi)\in C^\infty(\R^d\times \R^d)$ is an admissible cut-off function: there exists $\e_1,\e_2$ such that $0<\e_1<\e_2$ and
\beno
\chi(\theta,\eta)=1 \quad\textrm{for}\quad |\theta|\leq \e_1 |\eta|, \quad \chi(\theta,\eta)=0 \quad\textrm{for} \quad|\theta|\geq \e_2 |\eta|,
\eeno
and such that for any $(\theta,\eta)\in \R^d\times \R^d$,
\beno
|\pa_\theta^\al \pa_\eta^\beta\chi(\theta,\eta)|\leq C_{\al,\beta}(1+|\eta|)^{-|\al|-|\beta|}.
\eeno
The cut-off function $\psi(\eta)\in C^\infty(\R^d)$ satisfies
\beno
\psi(\eta)=0 \quad\textrm{for}\quad |\eta|\leq 1, \quad \psi(\eta)=1 \quad\textrm{for}\quad |\eta|\geq 2.
\eeno
Throughout this paper, we will take the admissible cut-off function $\chi(\theta,\eta)$ as
\beno
\chi(\theta,\eta)=\sum_{k=0}^\infty\zeta_{k-3}(\theta)\varphi_k(\eta),
\eeno
where $\zeta(\theta)=1$ for $|\theta|\le 1.1$ and $\zeta(\theta)=0$ for $|\theta|\ge 1.9$; and
\beno
&&\zeta_k(\theta)=\zeta(2^{-k}\theta)\quad \textrm{for}\quad k\in\Z,\\
&&\varphi_0=\zeta,\quad \varphi_k=\zeta_k-\zeta_{k-1}\quad \textrm{for}\quad k\ge 1.
\eeno

We also introduce the Littlewood-Paley operators $\Delta_k, S_k$ defined by
\beno
&&\Delta_k u={\mathcal F}^{-1}\big(\varphi_k(\xi)\widehat u(\xi)\big)\quad\textrm{ for }k\ge 0, \quad \Delta_ku=0\quad\textrm{ for }k<0,\\
&&S_ku={\mathcal F}^{-1}\big(\zeta_k(\xi)\widehat u(\xi)\big)\quad\textrm{ for }k\in\Z.
\eeno
In the case when the function $a$ depends only on the first variable $x$ in $T_au$, we  take $\psi=1$. Then $T_au$ is just the usual
Bony's paraproducts and
\ben\label{paraproduct}
T_au=\sum_{k=0}^{+\infty}S_{k-3}a\Delta_ku.
\een
Furthermore, we have Bony's decomposition:
\ben\label{Bony}
au=T_au+T_ua+R(u,a),
\een
where the remainder term  $R(u,a)$ is defined by
\beno
&&R(u,a)=\sum_{|k-\ell|\le 2; k,\ell\ge -2}\Delta_{k}a\Delta_\ell u.
\eeno

Now we introduce the Besov space.

\begin{definition}
Let $s\in\R$, $p,q\in [1,\infty]$. The inhomogeneous Besov space $B^s_{p,q}(\R^d)$
consists of the temperate distribution $f$ satisfying
\beno
\|f\|_{B^s_{p,q}}\eqdef \Big(\sum_{k}2^{ksq}\|\Delta_kf\|_{L^p}^q\Big)^\f1q<+\infty.
\eeno
\end{definition}

When $p=q=2$, $B^{s}_{p,q}(\R^d)$ is just the usual Sobolev space $H^s(\R^d)$; When $p=q=\infty$ and $s\notin \mathbf{N}$,
it is the Zygmund space $C^s(\R^d)$. From the definition, it is easy to see that if $s_1\le s_2$ and $q_2\le q_1$, then
\ben\label{Besov-p1}
\|f\|_{B^{s_1}_{p,q_1}}\le \|f\|_{B^{s_2}_{p,q_2}};
\een
if $s_2>s_1$, then
\ben\label{Besov-p2}
\|f\|_{B^{s_1}_{p,1}}\le \|f\|_{B^{s_2}_{p,\infty}}.
\een

The following Berstein's inequality will be repeatedly used.

\begin{lemma}\label{lem:Berstein}
Let $1\le p\le q\le \infty, \al\in \mathbf{N}^d$. Then it holds that
\beno
&&\|\pa^\al S_k u\|_{L^q}\le C2^{kd(\f 1p-\f1q+|\al|)}\|S_ku\|_{L^p}\quad\textrm{ for }k\in\mathbf{N},\\
&&\|\Delta_ku\|_{L^q}\le C2^{kd(\f 1p-\f1q-|\al|)}\sup_{|\beta|=|\al|}\|\pa^\beta\Delta_ku\|_{L^p}\quad\textrm{ for }k\ge 1.
\eeno

\end{lemma}

\subsection{Symbolic calculus}

We recall the symbolic calculus for the paradifferential operators.

\begin{proposition}\label{prop:symbolic calculus}
Let $m, m'\in \R$ and $\rho\in[0,1]$.

(1) If $a\in\Gamma^m_0(\mathbf{R}^d)$, then for any $\mu\in\R$,
\beno
\|T_a\|_{H^\mu\rightarrow H^{\mu-m}}\leq CM_0^m(a);
\eeno

(2) If $a\in\Gamma^m_\rho(\mathbf{R}^d), b \in\Gamma^{m'}_\rho(\mathbf{R}^d)$, then for any $\mu\in\R$,
\beno
\|T_aT_b-T_{ab}\|_{H^\mu\rightarrow H^{\mu-m-m'+\rho}}\leq CM_\rho^m(a)M_0^{m'}(b)+ K M_0^m(a)M_\rho^{m'}(b);
\eeno

(3) If $a\in\Gamma^m_\rho(\mathbf{R}^d)$, then for any $\mu\in\R$,
\beno
\|T_a-(T_a)^*\|_{H^\mu\rightarrow H^{\mu-m+\rho}}\leq CM_\rho^m(a).
\eeno
Here $(T_a)^*$ is the adjoint operator of $T_a$, and $C$ is a constant independent of $a,b$.
\end{proposition}

\begin{remark}\label{rem:SC-Holder}
If $\mu,\mu+m\not\in\mathbf{N}$, then we have
\beno
\|T_a\|_{W^{\mu+m,\infty}\rightarrow W^{\mu,\infty}}\leq CM_\rho^m(a),
\eeno
and if $\mu, \mu-m-m'+\rho\notin \mathbf{N}$, the we have
\beno
\|T_aT_b-T_{ab}\|_{W^{\mu,\infty}\rightarrow W^{\mu-m-m'+\rho,\infty} }\leq C M_\rho^m(a)M_0^{m'}(b)+ C M_0^m(a)M_\rho^{m'}(b).
\eeno
\end{remark}

\begin{lemma}\label{lem:Besov}
Let $m, m', \mu\in \R, q\in [1,\infty]$ and $\rho\in[0,1]$.

 (1)If $a\in\Gamma^m_0(\mathbf{R}^d)$, then
\beno
\|T_a\|_{{B^\mu_{\infty,q}}\rightarrow B^{\mu-m}_{\infty,q}}\leq C M_0^m(a);
\eeno

(2) If $a\in\Gamma^m_\rho(\mathbf{R}^d), b \in\Gamma^{m'}_\rho(\mathbf{R}^d)$, then
\beno
\|T_aT_b-T_{ab}\|_{B^\mu_{\infty,q}\rightarrow B^{\mu-m-m'+\rho}_{\infty,q}}\leq CM_\rho^m(a)M_0^{m'}(b)+ CM_0^m(a)M_\rho^{m'}(b).
\eeno
Here $C$ is a constant independent of $a,b$.
\end{lemma}

\no{\bf Proof.}\,Take $\al$ such that $\al,\al+m\notin \mathbf{N}$. From the definition of $T_a$, we know that there exists a constant $N_0\in \mathbf{N}$ such that
\ben\label{eq:localization}
\Delta_j(T_a u)=\sum_{|j-k|\leq N_0}\Delta_j(T_a \Delta_k u).
\een
Then it follows from Remark \ref{rem:SC-Holder} that
\begin{align*}
\|\Delta_j T_a u\|_{L^\infty}&\leq \sum_{|j-k|\leq N_0} \|\Delta_j T_a \Delta_k u\|_{L^\infty} \\
&\leq \sum_{|j-k|\leq N_0} 2^{-j\al}  \|T_a \Delta_k u\|_{W^{\al,\infty}}\leq  C\sum_{|j-k|\leq N_0} 2^{-j\al}  \| \Delta_k u\|_{W^{\al+m,\infty}}\\
&\leq  C \sum_{|j-k|\leq N_0}2^{mk}  \|\Delta_k u\|_{L^\infty},
\end{align*}
which implies (1).  The proof of (2) is similar.\ef

\begin{remark}\label{rem:symb}
If the symbol $a(x,\xi)$ satisfies
\beno
\|\pa_\xi^\al a(\cdot,\xi)\|_{C^{-\rho}}\leq C_\al(1+|\xi|)^{m-|\al|}\quad \textrm{for any} \quad |\xi|\geq \f12,
\eeno
for some $\rho>0$, then $T_a$ is bounded from $H^\mu(\R^d)$ to $H^{\mu-m-\rho}(\R^d)$ and $B^\mu_{\infty,q}$ to $B^{\mu-m-\rho}_{\infty,q}(\R^d)$.
\end{remark}

\begin{lemma}\label{lem:nonlinear}
Let $F$ be a smooth function with $F(0)=0$ and $\mu>0, q\in [1,\infty]$. Then it holds that
\beno
&&\|F(u)\|_{H^\mu}\le C(\|u\|_{L^\infty})\|u\|_{H^\mu},\\
&&\|F(u)\|_{B^\mu_{\infty,q}}\le C(\|u\|_{L^\infty})\|u\|_{B^\mu_{\infty,q}}.
\eeno
\end{lemma}

\subsection{Commutator estimates}

\begin{proposition}\label{prop:commu-Ds}
Let $m, \mu\in \R, s>0$ and $a\in \Gamma_\rho^m(\R^d)$ with $\rho\in (0,1]$. Then
\beno
\big\|[\Ds, T_a]u\big\|_{H^\mu}\le CM_\rho^m(a)\|u\|_{H^{s+\mu+m-\rho}}.
\eeno
\end{proposition}

\no{\bf Proof.}\,We write
\beno
[\Ds, T_a]u=T_{\langle\xi\rangle^s}T_au-T_{\langle\xi\rangle^sa}+(\Ds-T_{\langle\xi\rangle^s})T_au,
\eeno
then the proposition follows from Proposition \ref{prop:symbolic calculus} and
\beno
\|(\Ds-T_{\langle\xi\rangle^s})T_au\|_{H^\mu}\le C\|\Ds(1-\psi(D))T_au\|_{H^{\mu}}\le CM_0^m(a)\|u\|_{H^{-\mu'}}
\eeno
for any $\mu'>0$.\ef

\begin{proposition}\label{prop:commutator}
Let $V\in C([0,T]; B^1_{\infty,1}(\mathbf{R}^d))$ and $p=p(t,x,\xi)$ be homogenous in $\xi$ of order $m$. Then it holds that
\beno
\big\|[T_p,\pa_t+V\cdot \nabla ]u(t)\big\|_{L^2}
\leq C\big(M^m_0(p)\|V(t)\|_{B^1_{\infty,1}}+M^m_0(\pa_t p+V\cdot\nabla p)\big)\|u(t)\|_{H^m}.
\eeno
\end{proposition}

\no{\bf Proof.}\,We follow the proof of Lemma 2.16 in \cite{ABZ3}.
As in \cite{ABZ3}, it suffices to consider the case when $p=p(t, x)$ by decomposing $p$ into a sum of spherical harmonic.
By a direct calculation, we have
\begin{align*}
[\pa_t +V\cdot\nabla,T_p]u
&=T_{\pa_t p}u+V\cdot T_{\nabla p}u+V\cdot T_p\nabla u-T_p(V\cdot\nabla u)\\
&=T_{\pa_t p+V\cdot\nabla p}u+\big(V\cdot T_{\nabla p}-T_{V\cdot\nabla p}\big)u+\big(V\cdot T_p\nabla u-T_p(V\cdot\nabla u)\big).
\end{align*}

First of all, we get by Proposition \ref{prop:symbolic calculus} that
\beno
\|T_{\pa_t p+V\cdot\nabla p}u\|_{L^2}\le C\|\pa_t p+V\cdot\nabla p\|_{L^\infty}\|u\|_{L^2}.
\eeno
Set $S^{j-3}(V)=\sum_{k\geq j-3}\Delta_k V$. Then $V=S_{j-3}(V)+S^{j-3}(V)$ .
Hence,
\begin{align*}
V\cdot T_{\nabla p}u
=&\sum_j S_{j-3}(V)\cdot S_{j-3}(\nabla p)\Delta_j u+\sum_j S^{j-3}(V)\cdot S_{j-3}(\nabla p)\Delta_j u\\
=&\sum_{j}(S_{j-3}(S_{j-3}(V)\cdot\nabla p)\Delta_j u+([S_{j-3}(V),S_{j-3}]\cdot \nabla p)\Delta_j u)\\
&\quad+\sum_j S^{j-3}(V)\cdot S_{j-3}(\nabla p)\Delta_j u\\
=&T_{V\cdot\nabla p}u-\sum_{j}\big(S_{j-3}(S^{j-3}(V)\cdot\nabla p\big)+\big([S_{j-3}(V),S_{j-3}]\cdot \nabla p)\big)\Delta_j u\\
    & +\sum_j S^{j-3}(V)\cdot S_{j-3}(\nabla p)\Delta_j u.
\end{align*}
By Lemma \ref{lem:Berstein}, we get
\begin{align*}
\big\|\sum_j S^{j-3}(V)\cdot S_{j-3}(\nabla p)\Delta_j u\big\|_{L^2}
&\leq  C\|p\|_{L^\infty}\sum_j 2^j\|S^{j-3}(V)\|_{L^\infty}\|\Delta_j u\|_{L^2}\\
&\le C\|p\|_{L^\infty}\|u\|_{L^2}\sum_j2^j\sum_{k\ge j-3}\|\Delta_k V\|_{L^\infty}\\
&\leq C \|p\|_{L^\infty} \|V\|_{B^{1}_{\infty,1}}\|u\|_{L^2}.
\end{align*}
Noting that $S_{j-3}(S^{j-3}(V)\cdot\nabla p)\Delta_j u$ is spectrally supported in an annulus $\{c_12^j\leq|\xi|\leq c_2 2^j\}$,
we infer from Lemma \ref{lem:Berstein} that
\beno
&&\big\|\sum_{j}S_{j-3}(S^{j-3}(V)\cdot\nabla p)\Delta_j u\big\|_{L^2}^2\\
&&\leq \sum_j \big( 2^j\|S^{j-3}(V)p\|_{L^\infty}+\|S^{j-3}(\nabla \cdot V)p\|_{L^\infty}\big)^2\|\Delta_j u\|_{L^2}^2\\
&&\leq C\|p\|_{L^\infty}^2\|V\|_{B^{1}_{\infty,1}}^2\sum_j\|\Delta_j u\|_{L^2}\\
&&\leq C \|p\|_{L^\infty}^2\|V\|_{B^{1}_{\infty,1}}^2\|u\|_{L^2}^2,
\eeno
Similarly, $([S_{j-3}(V),S_{j-3}]\nabla p)\Delta_j u$ is also spectrally supported in an annulus $\{c_12^j\leq|\xi|\leq c_2 2^j\}$, thus,
\begin{align*}
\big\|\sum_{j}([S_{j-3}(V),S_{j-3}]\nabla p)\Delta_j u\big\|^2_{L^2}
&\leq C\sum_{j}\|([S_{j-3}(V),S_{j-3}]\nabla p)\Delta_j u\|^2_{L^2}\\
&\leq C\sum_j\|S_{j-3}(V)\|^2_{W^{1,\infty}}\|p\|_{L^\infty}^2\|\Delta_j u\|_{L^2}^2\\
&\leq C \|V\|^2_{W^{1,\infty}}\|p\|_{L^\infty}^2\|u\|_{L^2}^2.
\end{align*}
Here we used the commutator estimate
\beno
\|[S_j,g]\na f\|_{L^\infty}\le C\|\na g\|_{L^\infty}\|f\|_{L^\infty},
\eeno
which follows from the identity
\begin{align*}
[S_j,g]\na f(x)=&\int_{\R^d}\breve{\zeta}_k(x-x')(g(x')-g(x))\na f(x')dx'\\
=&\int_{\R^d}\na\breve{\zeta}_k(x-x')(g(x')-g(x))f(x')dx'\\
&-\int_{\R^d}\breve{\zeta}_k(x-x')\na g(x')f(x')dx',
\end{align*}
and $\|\breve{\zeta}_k\|_{L^1}+\|x\na\breve{\zeta}_k\|_{L^1}\le C$.
This proves that
\beno
\|\big(V\cdot T_{\nabla p}-T_{V\cdot\nabla p}\big)u\|_{L^2}\le C \|p\|_{L^\infty} \|V\|_{B^{1}_{\infty,1}}\|u\|_{L^2} .
\eeno

Next we write
\beno
V\cdot (T_p \nabla u)=\sum_j S_{j-3}(V)S_{j-3}(p)\cdot\Delta_j \nabla u+\sum_j S^{j-3}(V)S_{j-3}(p)\cdot\Delta_j \nabla u.
\eeno
It follows from Lemma \ref{lem:Berstein}  that
\beno
\big\|\sum_j S^{j-3}(V)S_{j-3}(p)\Delta_j \nabla u\big\|_{L^2}\le C \|p\|_{L^\infty} \|V\|_{B^{1}_{\infty,1}}\|u\|_{L^2} .
\eeno
On the other hand,  we have
\beno
&&T_p(V\cdot\nabla u)=T_pT_V\cdot \nabla u+T_p(V-T_V)\cdot\nabla u\\
&&=\sum_j S_{j-3}(p)\Delta_j\{S_{j-3}(V)\cdot \nabla u\}\\
&&\quad+\sum_{j,k} S_{j-3}(p)\Delta_j\{(S_{k-3}-S_{j-3})(V)\cdot\nabla \Delta_k u\}+T_p(V-T_V)\nabla u\\
&&= \sum_j S_{j-3}(V)S_{j-3}(p)\Delta_j \nabla u +\sum_jS_{j-3}(p)[\Delta_j,S_{j-3}(V)]\cdot\nabla u\\
&&\quad+\sum_{j,k} S_{j-3}(p)\cdot\Delta_j\{(S_{k-3}-S_{j-3})(V)\cdot\nabla \Delta_k u\}+T_p(V-T_V)\cdot\nabla u\\
&&\equiv \sum_j S_{j-3}(V)S_{j-3}(p)\cdot \na \Delta_j u+I_1+I_2+I_3.
\eeno
We get by Proposition \ref{prop:symbolic calculus} and (\ref{Bony}) that
\beno
\|I_3\|_{L^2}\le C \|p\|_{L^\infty}\|(V-T_V)\cdot\na u\|_{L^2}\le C \|p\|_{L^\infty} \|V\|_{B^{1}_{\infty,1}}\|u\|_{L^2} .
\eeno
Note that the summation index $(j,k)$ in $I_2$ should satisfy $|k-j|\le N_0$ for some $N_0\in \mathbf{N}$, hence,
\beno
\|I_2\|_{L^2}\le C \|p\|_{L^\infty} \|V\|_{B^{1}_{\infty,1}}\|u\|_{L^2} .
\eeno
We rewrite $I_1$ as
\beno
I_1=\sum_jS_{j-3}(p)[\Delta_j,S_{j-N_0}(V)]\cdot\nabla u+\sum_jS_{j-3}(p)[\Delta_j,(S_{j-3}-S_{j-N_0})(V)]\cdot\nabla u
\eeno
for some $N_0\in \mathbf{N}$ so that $S_{j-3}(p)[\Delta_j,S_{j-N_0}(V)]\cdot\nabla u$ is  spectrally supported in an annulus $\{c_12^j\leq|\xi|\leq c_2 2^j\}$.
Then as in the above, it is easy to get
\beno
\|I_1\|_{L^2}\le C \|p\|_{L^\infty} \|V\|_{B^{1}_{\infty,1}}\|u\|_{L^2} .
\eeno
Hence,  we conclude
\beno
\|V\cdot T_p\nabla u-T_p(V\cdot\nabla u)\|_{L^2}\le C \|p\|_{L^\infty} \|V\|_{B^{1}_{\infty,1}}\|u\|_{L^2} .
\eeno
This finishes the proof. \ef

\section{Parabolic evolution equation}

Let $I=[z_0, z_1]$. We denote by $\G_\r^m(I\times \R^d)$ the space of symbols $a(z;x,\xi)$ satisfying
\beno
\widetilde{M}^m_\rho(a)\eqdef\sup_{z\in I}\sup_{|\alpha|\leq 3d/2+1+\rho}\sup_{|\xi|\geq1/2} \|(1+|\xi|)^{|\alpha|-m}\pa^\alpha_\xi a(z;\cdot,\xi)\|_{W^{\rho,\infty}}
<+\infty.
\eeno

In this section, we study the parabolic evolution equation
\ben\label{eq:parabolic evolution}
\left\{
\begin{array}{l}
\pa_z w+T_p w=f,\\
{w|_{z=z_0}=w_0,}
\end{array}
\right.
\een
where the symbol $p\in \G^1_\r(I\times \R^d)$ is elliptic in the sense that there exists $c_1>0$ such that for any $z\in I, (x,\xi)\in \R^d\times \R^d$,
it holds that
\ben
\textrm{Re}\,p(z;x,\xi)\geq c_1|\xi|.
\een

In order to obtain the maximal parabolic regularity of the solution,
we introduce Chemin-Lerner type space $\widetilde{L}_z^q(I; B^r_{p,\ell}(\R^d))$, whose norm is defined by
\beno
\|f\|_{\widetilde{L}_z^q(I; B^r_{p,\ell})}\eqdef\Big(\sum_{k}2^{kr\ell}\|\Delta_k f\|_{L_z^q(I; L^p)}^\ell\Big)^\f1\ell,
\eeno
which was firstly introduced by Chemin and Lerner \cite{Che} to study the incompressible Navier-Stokes equations.
When $p=\ell=\infty$, we denote it by $\widetilde{L}_z^q(I; C^s(\R^d))$. When $p=q=\ell=2$,
we have $\widetilde{L}_z^q(I; B^r_{p,\ell}(\R^d))\equiv L^2(I;H^r(\R^d))$; When $q=\infty, p=\ell=2$, we denote it by $\widetilde{L}_z^\infty(I; H^r(\R^d))$. In this case, we have
\ben\label{chemin space}
\|f\|_{{L}_z^\infty(I; H^r)}\le \|f\|_{\widetilde{L}_z^\infty(I; H^r)}.
\een

\begin{proposition}\label{prop:parabolic estimate-maximal}
Let $r\in \R, \ell\in [1,\infty]$ $1\le q\le p\le \infty$.
Assume that $p\in \Gamma^1_\rho(I\times \R^d)$ for $\rho>0$ and $w$ is a solution of (\ref{eq:parabolic evolution}).
Then for any $\delta>0$, we have
\beno
\|w\|_{\widetilde{L}_z^p(I;B^{r+\f1p}_{\infty,\ell})}\leq C\big(\widetilde{M}^1_{\rho}(p), c_1^{-1}\big)
\Big(\|w_0\|_{B^r_{\infty,\ell}}+\|f\|_{\widetilde{L}_z^{q}(I;B^{r-1+\f1{q}}_{\infty,\ell})}+\|w\|_{\widetilde{L}_z^p(I;C^{-\delta})}\Big),
\eeno
where $C(\cdot)$ is a nondecreasing function independent of $p$.
\end{proposition}

The proof is based on the following classical parabolic smoothing effect.

\begin{Lemma}\label{lem:parabolic est-classical}
Let $\kappa>0$ and $p\in [1,\infty]$. Then there exists some $c>0$ such that for any $t>0, k\ge 1$, we have
\beno
\|e^{-t\kappa|D|}\Delta_ku\|_{L^p}\le Ce^{-ct2^{k}}\|\Delta_ku\|_{L^p}.
\eeno
\end{Lemma}

\no{\bf Proof.} Take a function $\chi_1\in C_0^\infty(\R^d\setminus \{0\})$ such that $\chi_1(\xi/2^k)=1$ for $\xi\in\textrm{supp}\varphi_k$.
Then we have
$$
e^{-t\kappa|D|}\Delta_ku(x)=(2\pi)^{-d}\int_{\R^d}e^{ix\cdot\xi}e^{-\kappa
t |\xi|}\chi_1(\xi/2^k)\widehat{\Delta_ku}(\xi)d\xi=G_{k,t}\ast \Delta_ku(x),
$$
where
$$
G_{k,t}(x)=(2\pi)^{-d}\int_{\R^d}e^{ix\cdot\xi}e^{-\kappa t
|\xi|}\chi_1(\xi/2^k)d\xi.
$$
Then the lemma will follows from  Young's inequality and the estimate
\beq\label{eq:kernel}
|G_{k,t}(x)|\le Ce^{-ct2^{k}}2^{dk}\big(1+|2^kx|\big)^{-N}
\eeq
for some $N>d$.

Now we prove (\ref{eq:kernel}). Noting that
$$
G_{k,t}(x)=(2\pi)^{-d}2^{dk}\int_{\R^d}e^{i2^kx\cdot\xi}e^{-\kappa
t2^{k} |\xi|}\chi_1(\xi)d\xi\equiv
2^{dk}\widetilde{G}_{k,t}(2^kx),$$
thus it suffices to show that
\beq\label{eq:kernel-est1}
\widetilde{G}_{k,t}(x)\le Ce^{-ct2^{k}}(1+|x|)^{-N}. \eeq
It is easy to see that
\begin{eqnarray}\label{eq:kernel-est2}
|\widetilde{G}_{k,t}(x)|\le
C\int_{|\xi|\thicksim 1}e^{-\kappa t2^{k}}d\xi\le Ce^{-ct2^{k}}.
\end{eqnarray}
To obtain the behavior of $\widetilde{G}_{k,t}(x)$ for large $x$,
we need to integrate by parts. For this end, we introduce the operator $L(x,D)=\frac {x\cdot \nabla_\xi} {i|x|^2}$.
Since $L(x,D)e^{ix\cdot \xi}=e^{ix\cdot\xi}$, then for any $N\in \mathbf{N}$, we have
\beno
\widetilde{G}_{k,t}(x)&=&(2\pi)^{-d}\int_{\R^d}L^N(e^{ix\cdot\xi})e^{-\kappa
t2^{k}|\xi|}\chi_1(\xi)d\xi\nonumber\\
&=&(2\pi)^{-d}\int_{\R^d}e^{ix\cdot\xi}(L^*)^N(e^{-\kappa
t2^{k} |\xi|}\chi_1(\xi))d\xi,
\eeno
where the integrand can be majorized by
$$
|x|^{-N}\max\big(1,(t2^{k})^N,t2^{k}|\xi|^{1-N}\big)e^{-\kappa
t2^{k}|\xi|}1_{|\xi|\thicksim 1}.$$
Hence, we infer that
\begin{align*}
|\widetilde{G}_{k,t}(x)|&\le C_N|x|^{-N}e^{-ct2^{k}}\int_{|\xi|\sim 1}(1+t2^{k}|\xi|)^Ne^{-\frac
{\kappa t2^{k}} 2|\xi|}d\xi\nonumber\\
&\le C_N|x|^{-N}e^{-ct2^{k}},
\end{align*}
which along with (\ref{eq:kernel-est2}) implies (\ref{eq:kernel-est1}).\ef

\vspace{0.2cm}

\no{\bf Proof of Proposition \ref{prop:parabolic estimate-maximal}.}\,For $y\in I$ and $z\in [z_0,y]$, we set
\beno
e(y,z;x,\xi)=\exp\Big(-\int_z^y p(s;x,\xi) ds\Big).
\eeno
Noting that $\pa_ze=ep$, we get by (\ref{eq:parabolic evolution}) that
\beno
\pa_z(T_ew)=T_{\pa_ze}w+T_e\pa_zw=\big(T_{ep}-T_eT_p\big)w+T_ef.
\eeno
Integrating it on $[z_0,y]$, we get
\begin{align*}
T_1w(y)&=T_{e|_{z=z_0}}w_0+\int^y_{z_0}T_ef(z)dz+\int^y_{z_0}(T_{ep}-T_eT_p)w(z)dz\\
&\triangleq G_1+G_2+Rw
\end{align*}
so that for any $N\in \mathbf{N}$, there holds
\ben\label{eq:w}
w=(I+R+\cdots+R^N)(G_1+G_2-T_1w+w)+R^{N+1}w,
\een
where for any $\delta>0$, we have
\beno
\|w-T_1w\|_{\widetilde{L}_y^p(I; B^{r+\f1p}_{\infty,\ell})}
=\|(1-\psi(D))w\|_{\widetilde{L}_y^p(I; B^{r+\f1p}_{\infty,\ell})} \leq C\|w\|_{\widetilde{L}_y^p(I; C^{-\delta})}.
\eeno

It is east to verify that $e(y,z;x,\xi)\exp({c_1}(y-z)|\xi|/2)\in \G^0_\rho(\R^d)$ for $y,z\in I, z\le y$ with the bound
\beno
M_\rho^0\big(e(y,z;x,\xi)\exp({c_1}(y-z)|\xi|/2)\big)\le C\widetilde{M}_\rho^1(p).
\eeno
Thus by (\ref{eq:localization}), Remark \ref{rem:SC-Holder} and Lemma \ref{lem:parabolic est-classical}, we have
\begin{align}\label{eq:parabolic est}
\|\Delta_j T_e u\|_{L^\infty} &\leq \sum_{|j-k|\leq N_0}\|\Delta_j T_{e}\Delta_k u\|_{L^\infty}\nonumber\\
&\leq \sum_{|j-k|\leq N_0}2^{-\f j2}\|\Delta_jT_{e}\exp(c_1|D|(y-z)/2)\exp(-c_1|D|(y-z)/2) \Delta_k u \|_{C^{\f12}}\nonumber\\
&\leq C\widetilde{M}_\rho^1(p)\sum_{|j-k|\leq N_0} 2^{-\f j2}\|\exp(-c_1|D|(y-z)/2) \Delta_k u \|_{C^{\f12}}\nonumber\\
&\leq C\widetilde{M}_\rho^1(p)\sum_{|j-k|\leq N_0}\exp(-c(y-z)2^k)\|\Delta_k u\|_{L^\infty}
\end{align}
for some $N_0\in \mathbf{N}$ and $c>0$(\textbf{Important note}: the summation index $k\ge 1$ due to the definition of $T_e$).

Now let us turn to the estimates of $G_i$. We get by (\ref{eq:parabolic est}) that
\begin{align}
&\|G_1\|_{\widetilde{L}_y^p(I; B^{r+\f1p}_{\infty,\ell})}
=\Big(\sum_j 2^{j\ell(r+\f1p)}\|\Delta_j G_1\|_{L_y^p(I; L^\infty)}^\ell\Big)^\f1\ell\nonumber\\
&\leq C\widetilde{M}_\rho^1(p)\Big(\sum_j\big(\sum_{|j-k|\leq N_0}2^{j(r+\f1p)}\big\|\exp(-c(y-z_0)2^k)\|\Delta_kw_0\|_{L_z^\infty}\big\|_{L_y^p(I)}\big)^\ell\Big)^\f1\ell\nonumber \\
&\le  C\widetilde{M}_\rho^1(p)\Big(\sum_j\big(\sum_{|j-k|\leq N_0}2^{j(r+\f1p)}2^{-\f kp}\|\Delta_kw_0\|_{L_z^\infty}\big)^\ell\Big)^\f1\ell\nonumber \\
&\leq C\widetilde{M}_\rho^1(p)\|w_0\|_{B^r_{\infty,\ell}}.\label{eq:G1}
\end{align}
For $G_2$, we have by (\ref{eq:parabolic est}) that
\beno
 \int^y_{z_0}\|\Delta_j T_ef(z)\|_{L^\infty}dz\le C\widetilde{M}_\rho^1(p)\sum_{|j-k|\leq N_0}\int_{z_0}^y\exp(-c(y-z)2^k)\|\Delta_k f(z)\|_{L^\infty}dz,
\eeno
from which and Young's inequality, we infer that
\beno
\Big\|\int^y_{z_0}\|\Delta_j T_e f(z)\|_{L^\infty}dz\Big\|_{L^{p}(I)}
\leq C\widetilde{M}_\rho^1(p)\sum_{|j-k|\leq N_0}2^{-k(1+\f1p-\f1q)}\|\Delta_k f\|_{L^q_z(I;L^\infty)}.
\eeno
This implies that
\ben\label{eq:G2}
\|G_2\|_{\widetilde{L}_y^p(I;B^{r+\f1p}_{\infty,\ell})}\leq C\widetilde{M}_\rho^1(p)\|f\|_{\widetilde{L}_z^q(I; B^{r-1+\f1q}_{\infty,\ell})}.
\een

Similar to the proof of (\ref{eq:parabolic est}), we can get
\beno
\|\Delta_j(T_eT_p-T_{ep})w\|_{L^\infty}\leq C(\widetilde{M}^1_\rho(p))\sum_{|j-k|\leq N_0}2^{k(1-\rho)}\exp(-c(y-z)2^k)\|\Delta_k w\|_{L^\infty},
\eeno
which implies that
\ben\label{eq:Rw}
\|Rw\|_{\widetilde{L}_y^p(I;B^{r+\f1 p}_{\infty,\ell})}\leq C(\widetilde{M}^1_\rho(p))\|w\|_{\widetilde{L}_z^p(I;B^{r+\f1p-\rho}_{\infty,\ell})}.
\een

Take $N$ big enough so that $r+\f1p-(N+1)\rho\le -\delta$. Then the proposition follows from
(\ref{eq:w}) and (\ref{eq:G1})-(\ref{eq:Rw}).\ef

Given $r\in\R$,  let us introduce the spaces
\beno
&&X^r(I)\eqdef  \widetilde{L}^\infty_{z}(I;H^r(\mathbf{R}^d))\cap L^2_{z}(I;H^{r+\f12}(\mathbf{R}^d)),\\
&&Y^r(I)\eqdef \widetilde{L}^1_{z}(I;H^r(\mathbf{R}^d))+L^2_{z}(I;H^{r-\f12}(\mathbf{R}^d)).
\eeno

In a similar way as in Proposition \ref{prop:parabolic estimate-maximal}, one can show that

\begin{proposition}\label{prop:parabolic estimate}
Let $r\in \R$. Assume that $p\in \Gamma^1_\rho(I\times \R^d)$ for $\rho>0$ and $w$ is a solution of (\ref{eq:parabolic evolution}).
Then it holds that
\beno
\|w\|_{X^r(I)} \leq C\big(\widetilde{M}^1_{\rho}(p), c_1\big)\big(\|w_0\|_{H^r}+\|f\|_{Y^r(I)}+\|w\|_{L^2(I;H^r)}\big),
\eeno
where $C(\cdot)$ is a nondecreasing function independent of $p$.
\end{proposition}

Let us conclude this section by presenting some product estimates in the Chemin-Lerner type space.

\begin{lemma}\label{lem:PP-Hs}
Let $r\in \R$ and $q,q_1,q_2\in [1,\infty]$ with $\f 1 q=\f1 {q_1}+\f1 {q_2}$. Then for any $r_1, r_2>0$, we have
\beno
&&\|T_gf\|_{\widetilde{L}^q_z(I;H^r)}\le C\|g\|_{\widetilde{L}^{q_1}_z(I;L^\infty)}\|f\|_{\widetilde{L}^{q_2}_z(I;H^r)},\\
&&\|T_gf\|_{\widetilde{L}^q_z(I;H^r)}\le C\|g\|_{\widetilde{L}^{q_1}_z(I;C^{-r_1})}\|f\|_{\widetilde{L}^{q_2}_z(I;H^{r+r_1})},\\
&&\|T_gf\|_{\widetilde{L}^q_z(I;H^r)}\le C\|g\|_{\widetilde{L}^{q_1}_z(I;C^{0})}\|f\|_{\widetilde{L}^{q_2}_z(I;H^{r+r_2})}.
\eeno
\end{lemma}

\no{\bf Proof.}\,By the definition of paraproduct, we have
\beno
\Delta_jT_{g}f=\sum_{|j-k|\le N_0}\Delta_j\big(S_{k-3}g\Delta_kf\big)\quad\textrm{ for some }N_0\in \mathbf{N}.
\eeno
Hence, we get by Lemma \ref{lem:Berstein} that
\begin{align*}
\|\Delta_jT_{g}f\|_{\widetilde{L}^q_z(I;L^2)}&\le C\sum_{|j-k|\le N_0}\|S_{k-3}g\|_{\widetilde{L}^{q_1}_z(I;L^\infty)}\|\Delta_kf\|_{\widetilde{L}^{q_2}_z(I;L^2)}\\
&\le C\sum_{|j-k|\le N_0}\|g\|_{\widetilde{L}^{q_1}_z(I;L^\infty)}\|\Delta_kf\|_{\widetilde{L}^{q_2}_z(I;L^2)},
\end{align*}
which implies the fisrt inequality of the lemma. On the other hand, by the definition of $S_k$, we have
\begin{align*}
&\|S_{k-3}g\|_{\widetilde{L}^{q_1}_z(I;L^\infty)}\le \sum_{\ell\le k-2}\|\Delta_\ell g\|_{\widetilde{L}^{q_1}_z(I;L^\infty)}
\le 2^{jr_1}\|g\|_{\widetilde{L}^{q_1}_z(I;C^{-r_1})}\quad \textrm{or}\\
&\|S_{k-3}g\|_{\widetilde{L}^{q_1}_z(I;L^\infty)}\le Ck\|g\|_{\widetilde{L}^{q_1}_z(I;C^{0})}\le C2^{kr_2}\|g\|_{\widetilde{L}^{q_1}_z(I;C^{0})},
\end{align*}
which imply the last two inequalities.\ef\vspace{0.1cm}

In a similar way, one can show that
\begin{lemma}\label{lem:PP-Holder}
Let $r\in \R$ and $q,q_1,q_2,\ell\in [1,\infty]$ with $\f 1 q=\f1 {q_1}+\f1 {q_2}$. Then for any $r_1, r_2>0$, we have
\beno
&&\|T_gf\|_{\widetilde{L}^q_z(I;B^r_{\infty,\ell})}\le C\|g\|_{\widetilde{L}^{q_1}_z(I;L^\infty)}\|f\|_{\widetilde{L}^{q_2}_z(I;B^r_{\infty,\ell})},\\
&&\|T_gf\|_{\widetilde{L}^q_z(I;B^r_{\infty,\ell})}\le C\|g\|_{\widetilde{L}^{q_1}_z(I;C^{-r_1})}\|f\|_{\widetilde{L}^{q_2}_z(I;B^{r+r_1}_{\infty,\ell})},\\
&&\|T_gf\|_{\widetilde{L}^q_z(I;B^r_{\infty,\ell})}\le C\|g\|_{\widetilde{L}^{q_1}_z(I;C^{0})}\|f\|_{\widetilde{L}^{q_2}_z(I;B^{r+r_2}_{\infty,\ell})}.
\eeno
\end{lemma}

\begin{lemma}\label{lem:PR}
Let $q,q_1,q_2,\ell\in [1,\infty]$ with $\f 1 q=\f1 {q_1}+\f1 {q_2}$. Then for any $r>0$ and $r_1\in R$, we have
\beno
&&\|R(f,g)\|_{\widetilde{L}^q_z(I;H^r)}\le C\|g\|_{\widetilde{L}^{q_1}_z(I;C^{r_1})}\|f\|_{\widetilde{L}^{q_2}_z(I;H^{r-r_1})},\\
&&\|R(f,g)\|_{\widetilde{L}^q_z(I;B^r_{\infty,\ell})}\le C\|g\|_{\widetilde{L}^{q_1}_z(I;C^{r_1})}\|f\|_{\widetilde{L}^{q_2}_z(I;B^{r-r_1}_{\infty,\ell})}.
\eeno
If $r\le 0$ and $r_1+r_2>0$, then we have
\beno
&&\|R(f,g)\|_{\widetilde{L}^q_z(I;H^r)}\le C\|g\|_{\widetilde{L}^{q_1}_z(I;C^{r_1})}\|f\|_{\widetilde{L}^{q_2}_z(I;H^{r_2})},\\
&&\|R(f,g)\|_{\widetilde{L}^q_z(I;B^r_{\infty,\ell})}\le C\|g\|_{\widetilde{L}^{q_1}_z(I;C^{r_1})}\|f\|_{\widetilde{L}^{q_2}_z(I;C^{r_2})}.
\eeno
\end{lemma}

\no{\bf Proof.}\,Due to the definition of $R(f,g)$, we have
\beno
\Delta_jR(f,g)=\sum_{|k-\ell|\le 2; k,\ell\ge j-N_0}\Delta_j\big(\Delta_kf\Delta_\ell g\big)\quad\textrm{ for some }N_0\in \mathbf{N},
\eeno
from which and Lemma \ref{lem:Berstein}, we infer that
\begin{align*}
\|\Delta_jR(f,g)\|_{\widetilde{L}^q_z(I;L^2)}&\le C\sum_{|k-\ell|\le 2; k,\ell\ge j-N_0}\|\Delta_kf\|_{\widetilde{L}^{q_2}_z(I;L^2)}\|\Delta_\ell g\|_{\widetilde{L}^{q_1}_z(I;L^\infty)}\\
&\le C\|g\|_{\widetilde{L}^{q_1}_z(I;C^{r_1})}\sum_{k\ge j-N_0}2^{-kr_1}\|\Delta_kf\|_{\widetilde{L}^{q_2}_z(I;L^2)},
\end{align*}
which implies the first inequality of the lemma. The proof of the other three inequalities are similar. \ef

\section{Elliptic estimates in a strip of infinite depth}

In this section, we consider the elliptic equation in a strip of infinite depth $\cS=\big\{(x,y): x\in \R^d, y<\eta(x)\big\}$:
\begin{equation}\label{eq:elliptic}
\left\{\begin{aligned}
&\triangle_{x,y} \phi=g\quad\textrm{ in }\quad \cS,\\
&\phi|_{y=\eta}=f.
\end{aligned}\right.
\end{equation}
Throughout this section, we assume that $\eta\in H^{s+\f12}(\R^d)\cap C^{\f32+\e}(\R^d)$ for $s>1+\f d2$ and some $\e>0$.
We denote by $K_\eta=K_\eta\big(\|\eta\|_{C^{\f32+\e}},\|\eta\|_{L^2}\big)$ a nondecreasing function,
which may be different from line to line, $I=(-\infty,0)$.

First of all, we flatten the boundary of $\cS$ by the following regularized mapping:
\beno
(x,z)\in \R^d\times (-\infty,0]\longmapsto (x,\rho_\delta(x,z))\in\cS,
\eeno
where $\rho_\delta$ with $\delta>0$ is given by
\begin{eqnarray}\label{mapping}
\rho_\delta(x,z)=z+(e^{\delta z |D|}\eta)(x).
\end{eqnarray}
\begin{remark}
For any $z<0$, we have
\begin{align*}
\|\pa_z\rho_\delta-1\|_{L^\infty}&\le \delta\|e^{\delta z |D|}|D|\eta\|_{L^\infty}
=\delta\big\|P_{-\delta z}\ast |D|\eta\big\|_{L^\infty}\\
&\le C\delta\|P_{-\delta z}(\cdot)\|_{L^1}\||D|\eta\|_{L^\infty}
\le C\delta\|\eta\|_{C^{\f32+\e}}.
\end{align*}
Here $P_z(x)$ is the poisson kernel. Throughout this paper, we will fix $\delta$ small enough depending only on $\|\eta\|_{C^{\f32+\e}}$ such that
\beno
\|\pa_z\rho_\delta-1\|_{L^\infty}\le \f12, \quad\textrm{hence}\quad\pa_z\rho_\delta\ge \f12.
\eeno
\end{remark}

We set $v(x,z)=\phi(x,\rho_\delta(x,z))$. It is easy to find that $v$ satisfies
\begin{equation}\label{eq:ellitic-flat}
\left\{\begin{aligned}
&\pa^2_z v+\al \triangle v+\beta\cdot\nabla \pa_z v-\gamma \pa_z v=F_0,\\
&v|_{z=0}=f,
\end{aligned}\right.
\end{equation}
where $F_0=\alpha g$ and the coefficients $\al,\beta,\gamma$ are defined by
\ben\label{eq:elliptic coefficients}
\alpha=\f{(\pa_z \rho_\delta)^2}{1+|\nabla \rho_\delta|^2},\quad \beta=-2\f{\pa_z\rho_\delta \nabla \rho_\delta}{1+|\nabla \rho_\delta|^2},
\quad \gamma=\f{1}{\pa_z\rho_\delta}(\pa_z^2\rho_\delta+\alpha\Delta\rho_\delta+\beta\cdot\nabla\pa_z\rho_\delta).
\een

By the definition of $\rho_\delta$, we find
\ben\left\{\begin{aligned}
&\pa_z(\rho_\delta-z)-\delta |D|(\rho_\delta-z)=0,\\
&\rho_\delta-z|_{z=0}=\eta.\nonumber
\end{aligned}\right.
\een
Then we infer from Proposition \ref{prop:parabolic estimate} and Proposition \ref {prop:parabolic estimate-maximal}  that
\ben
&&\|\nabla \rho_\delta\|_{X^{s-\f12}(I)}+\|\pa_z\rho_\delta-1\|_{X^{s-\f12}(I)}\leq C\big(\|\eta\|_{C^{\f32+\e}}\big)\|\eta\|_{H^{s+\f12}},\label{eq:rho-Hs}\\
&&\|\nabla_{x,z}\rho_\delta\|_{\widetilde{L}^\infty_z(I;C^{\f12+\e})}
+\|\pa_{z}^2\rho_\delta\|_{\widetilde{L}^\infty_z(I;C^{-\f12+\e})}\leq C\big(\|\eta\|_{C^{\f32+\e}}\big).\label{eq:rho-Holder}
\een

In order to obtain the tame elliptic estimates, we paralinearize the elliptic equation (\ref{eq:ellitic-flat}) as
\ben\label{eq:elliptic-pl}
\pa^2_z v+T_\al\triangle v+T_\beta\cdot\nabla \pa_z v=F_0+F_1+F_2,
\een
with $F_1, F_2$ given by
\beno
F_1=\gamma \pa_zv,\quad F_2=(T_\al-\al)\triangle v+(T_\beta-\beta)\cdot\nabla \pa_z v.
\eeno
As in \cite{ABZ3}, the equation (\ref{eq:elliptic-pl}) can be decoupled into a forward and a backward parabolic
evolution equations:
\begin{eqnarray}\label{eq:elliptic-decouple}
(\pa_z-T_a)(\pa_z-T_A)v=F_0+F_1+F_2+F_3\triangleq F,
\end{eqnarray}
where
\beno
&&a=\f12\big(-i\beta\cdot\xi-\sqrt{4\al|\xi|^2-(\beta\cdot\xi)^2}\big),\\
&&A=\f12\big(-i\beta\cdot\xi+\sqrt{4\al|\xi|^2-(\beta\cdot\xi)^2}\big),\\
&&F_3=(T_aT_A-T_\al\triangle)v-(T_a+T_A+T_\beta\cdot\nabla )\pa_zv-T_{\pa_z A}v.
\eeno
\begin{remark}\label{rem:symbol-aA}
The symbols $a, A$ satisfy
\beno
a(z;x,\xi)\cdot A(z;x,\xi)=-\al(x,z)|\xi|^2,\quad a(z;x,\xi)+ A(z;x,\xi)=-i\beta(x,z)\cdot \xi.
\eeno
Noticing that
\beno
4\al|\xi|^2-(\beta\cdot\xi)^2\ge c_2|\xi|^2
\eeno
for some $c_2>0$ depending only on $\|\eta\|_{C^{\f32+\e}}$, it follows from (\ref{eq:rho-Holder}) that
\beno
\widetilde{M}^1_{\f12+\e}(a)\le C\big(\|\eta\|_{C^{\f32+\e}}\big),\quad \widetilde{M}^1_{\f12+\e}(A)\le C\big(\|\eta\|_{C^{\f32+\e}}\big).
\eeno
\end{remark}

\subsection{Elliptic estimates in Sobolev space}

\begin{proposition}\label{prop:elliptic Hs est}
Let $v$ be a solution of (\ref{eq:ellitic-flat}) on $I\times \R^d$. Then for all $\sigma\in [-\f12, s-\f12]$, it holds that
\begin{eqnarray*}
\|\nabla_{x,z} v\|_{X^{\sigma}(I)}\leq K_\eta\big(\|\na_{x,z}v\|_{L^2(I\times \mathbf{R}^d)}+
\|f\|_{H^{\sigma+1}}+ \|F_0\|_{Y^{\sigma}(I)}+\|\eta\|_{H^{s+\f12}}\|\nabla_{x,z} v\|_{L^{\infty}(I\times\R^d)}\big).
\end{eqnarray*}
Moreover, for $\sigma=-\f12$, we have
\beno
\|\na_{x,z}v\|_{X^{-\f12}(I)}\le K_\eta\big(\|F_0\|_{Y^{-\f12}(I)}+\|\nabla_{x,z}v\|_{L^2(I\times\R^d)}\big).
\eeno
\end{proposition}

Before proving the proposition, we make the estimates for the coefficients $\al,\beta, \gamma$ and $F_i(i=1, 2, 3)$.

\begin{lemma}\label{lem:coeffi}
It holds that
\beno
&&\|\alpha -1\|_{X^{s-\f12}(I)}+\|\beta\|_{X^{s-\f12}(I)}+\|\gamma\|_{X^{s-\f32}(I)}
\leq K_\eta\|\eta\|_{H^{s+\f12}},\\
&&\|\alpha\|_{\widetilde{L}_z^\infty(I; C^{\f12+\e})}+\|\beta\|_{\widetilde{L}_z^\infty(I; C^{\f12+\e})}+\|\gamma\|_{\widetilde{L}_z^\infty(I; C^{-\f12+\e})}\le K_\eta,\\
&&\|\alpha-1\|_{\widetilde{L}^1_z(I; H^{s+\f12})}+\|\beta\|_{\widetilde{L}^1_z(I; H^{s+\f12})}+\|\gamma\|_{\widetilde{L}^1_z(I; H^{s-\f12})}\le K_\eta\|\eta\|_{H^{s+\f12}},\\
&&\|\alpha\|_{\widetilde{L}_z^2(I; C^{1+\e})}+\|\beta\|_{\widetilde{L}_z^2(I; C^{1+\e})}
+\|\gamma\|_{\widetilde{L}_z^2(I; C^{\e})}\leq K_\eta.
\eeno
\end{lemma}

\no{\bf Proof.}\,Noting $s-\f32>0$, the first two inequalities of the lemma follows from Lemma \ref{lem:nonlinear} and (\ref{eq:rho-Hs})-(\ref{eq:rho-Holder})
except that $\|\gamma\|_{L_z^\infty(I; C^{-\f12+\e})}$. Thanks to (\ref{eq:rho-Holder}) and Lemma \ref{lem:nonlinear}, $\gamma$ can be written as
\beno
\gamma=\gamma_1\na^2\gamma_2,\quad \textrm{where}\quad \gamma_1\in \widetilde{L}_z^\infty(I;C^{\f12+\e}),\, \gamma_2\in \widetilde{L}^\infty_z(I;C^{\f32+\e})
\eeno
with the following bounds
\beno
\|\gamma_1\|_{\widetilde{L}_z^\infty(I;C^{\f12+\e})}+\|\gamma_2\|_{\widetilde{L}^\infty_z(I;C^{\f32+\e})}\leq K_\eta.
\eeno
We use Bony's decomposition (\ref{Bony}) to write $\gamma_1\na^2\gamma_2$ as
\beno
\gamma=T_{\gamma_1}\na^2\gamma_2+T_{\na^2\gamma_2}\gamma_1+R(\gamma_1,\na^2\gamma_2).
\eeno
We infer from Lemma \ref{lem:PP-Holder} and Lemma \ref{lem:PR} that
\beno
&&\|T_{\gamma_1}\na^2\gamma_2\|_{\widetilde{L}_z^\infty(I; C^{-\f12+\e})}\le C\|\gamma_1\|_{L^\infty}\|\gamma_2\|_{\widetilde{L}_z^\infty(I; C^{{\f32}+\e})}\le K_\eta,\\
&&\|T_{\na^2\gamma_2}\gamma_1\|_{\widetilde{L}_z^\infty(I; C^{2\e})}\le\|\gamma_1\|_{\widetilde{L}^\infty_z(I;C^{\f12+\e})}\|\gamma_2\|_{\widetilde{L}_z^\infty(I; C^{\f32+\e})} \le K_\eta,\\
&&\|R(\gamma_1,\na^2\gamma_2)\|_{\widetilde{L}_z^\infty(I; C^{2\e})}\le\|\gamma_1\|_{\widetilde{L}^\infty_z(I;C^{\f12+\e})}\|\gamma_2\|_{\widetilde{L}_z^\infty(I; C^{\f32+\e})}\le K_\eta.
\eeno
This gives the estimate of $\|\gamma\|_{\widetilde{L}_z^\infty(I; C^{-\f12+\e})}$.

For the last two inequalities of the lemma,  the proof is similar, but we need to use the following estimates for $\rho_\delta$:
\ben
&&\|\na^2_{x,z}\rho_\delta\|_{\widetilde{L}^1_z(I; H^{s-\f12})}\le C\|\eta\|_{H^{s+\f12}},\\
&&\|\nabla_{x,z}\rho_\delta\|_{\widetilde{L}^2_z(I;C^{1+\e})}
+\|\pa_{z}^2\rho_\delta\|_{\widetilde{L}^2_z(I;C^{\e})} \leq K_\eta.\label{eq:rho-H-L2}
\een
Indeed, by the the definition of $\rho_\delta$, we have
\begin{align*}
\|\na^2_{x,z}\rho_\delta\|_{\widetilde{L}^1_z(I; H^{s-\f12})}^2&\le \sum_j2^{2j(s-\f12)}\|\Delta_j\na^2_{x,z}\rho_\delta\|_{L^1_z(I;L^2)}^2\\
&\le C\sum_j2^{2j(s-\f32)}\|\Delta_j |D|^2e^{\delta z|D|}\eta\|_{L^1_z(I;L^2)}^2\\
&\le C\sum_j2^{2j(s+\f12)}\|\Delta_j\eta\|_{L^2}^2\le C\|\eta\|_{H^{s+\f12}},
\end{align*}
and by Lemma \ref{lem:parabolic est-classical}, we get
\begin{align*}
\|\nabla_{x,z}(\rho_\delta-z)\|_{\widetilde{L}^2_z(I;C^{1+\e})}&\le \sup_j2^{j(1+\e)}\|\Delta_j\na_{x,z}\rho_\delta\|_{L^2_z(I;L^\infty)}\\
&\le  \sup_j2^{j(1+\e)}\|\Delta_j |D|e^{\delta z|D|}\eta\|_{L^2_z(I;L^\infty)}\\
&\le  \sup_{j>0}2^{j(2+\e)}\|e^{cz\delta 2^j}\|_{L^2(I)}\|\Delta_j\eta\|_{L^\infty}+\|\Delta_0|D|e^{\delta z|D|}\eta\|_{L^2_z(I;L^\infty)}\\
&\le C\|\eta\|_{C^{\f 32+\e}}+C\|\eta\|_{L^2},
\end{align*}
where we use the estimate in the last inequality:
\beno
\|\Delta_0|D|e^{\delta z|D|}\eta\|_{L^2_z(I;L^\infty)}
\le C\|e^{\delta z|D|}|D|\Delta_0\eta\|_{L^2(I\times \R^d)}
\leq C\|\eta\|_{L^2}.
\eeno
The proof is finished. \ef

\begin{lemma}\label{lem:F1}
For any $0<\delta_1\leq \f12$ and $\sigma+\delta_1\leq s-\f12$, it holds that
\beno
\|F_1\|_{Y^{\sigma+\delta_1}(I)}\leq K_\eta
\big(\|\pa_zv\|_{L_z^2(I; H^{\sigma+\delta_1-\e})}+\|\pa_z v\|_{L^{\infty}(I\times\R^d)}\|\eta\|_{H^{s+\f12}}\big).
\eeno
\end{lemma}

\no{\bf Proof.}\,Using Bony's decomposition (\ref{Bony}), we write $F_1$ as
\[
F_1=\gamma\pa_zv=T_\gamma \pa_z v+T_{\pa_z v} \gamma +R(\gamma,\pa_z v).
\]
We infer from Lemma \ref{lem:PP-Hs} that
\begin{align*}
&\|T_\gamma \pa_z v\|_{L^2_z(I; H^{\sigma+\delta_1-\f12})}\leq C\|\gamma\|_{L_z^\infty(I; C^{-\f12+\e})}\|\pa_z v\|_{L_z^2(I; H^{\sigma+\delta_1-\e})},\\
&\|T_{\pa_z v}\gamma\|_{L^2_z(I; H^{\sigma+\delta_1-\f12})}\leq C\|\pa_z v\|_{L^{\infty}(I\times\R^d)}\|\gamma\|_{L_z^2(I; H^{\sigma+\delta_1-\f12})},
\end{align*}
and by noting $\sigma+\delta_1\le s-\f12$,
\beno
&\|R(\gamma,\pa_z v)\|_{\widetilde{L}^1_z(I; H^{\sigma+\delta_1})}\le
C\|\pa_z v\|_{\widetilde{L}_z^\infty(I; C^0)}\|\gamma\|_{\widetilde{L}_z^1(I; H^{s-\f12})}.
\eeno
This together with Lemma \ref{lem:coeffi} gives the lemma. \ef

\begin{lemma}\label{lem:F2}
For any $0<\delta_1\leq \f12$ and $\sigma+\delta_1\leq s-\f12$, it holds that
\beno
\|F_2\|_{Y^{\sigma+\delta_1}(I)}\leq K_\eta\|\nabla_{x,z}v\|_{L_z^\infty(I\times\R^d)}\|\eta\|_{H^{s+\f12}}.
\eeno

\end{lemma}

\no{\bf Proof.}\,Recalling  $F_2=(T_\al-\al)\triangle v+(T_\beta-\beta)\cdot\nabla \pa_z v$, it suffices to consider $(T_\al-\al)\triangle v$.
We get by (\ref{Bony}) that
\beno
(T_\al-\al)\triangle v=-T_{\Delta v}(\al-1)-R(\Delta v,\al-1).
\eeno
Due to $\sigma+\delta_1\le s-\f12$, we infer from Lemma \ref{lem:PP-Holder} and Lemma \ref{lem:PR} that
\begin{align*}
\|T_{\triangle v}(\al-1)\|_{L_z^2(I; H^{\sigma+\delta_1-\f12})}
&\leq C\|\nabla v\|_{L_z^\infty(I\times\R^d)} \|\al-1\|_{L^2_z(I; H^{\sigma+\delta_1+\f12})}\\
&\leq C\|\nabla v\|_{L_z^\infty(I\times\R^d)} \|\al-1\|_{L^2_z(I; H^{s})},\\
\|R(\Delta v,\al-1)\|_{\widetilde{L}_z^1(I; H^{\sigma+\delta_1})}&\le \|R(\Delta v,\al-1)\|_{\widetilde{L}_z^1(I; H^{s-\f12})}\\
&\le C\|\nabla v\|_{\widetilde{L}_z^\infty(I; C^0)}\|\al-1\|_{\widetilde{L}^1_z(I; H^{s+\f12})},
\end{align*}
which along with Lemma \ref{lem:coeffi} give the lemma.\ef

\begin{lemma}\label{lem:F3}
For any $\sigma,\delta_1\in \R$, it holds that
\beno
\|F_3\|_{Y^{\sigma+\delta_1}(I)}\leq K_\eta\|\nabla v\|_{L_z^2(I; H^{\sigma+\delta_1-\e})}.
\eeno
\end{lemma}

\no{\bf Proof.}\,It follows from Remark \ref{rem:symbol-aA} and Proposition \ref{prop:symbolic calculus} that
\begin{eqnarray*}
&&\|(T_aT_A-T_\al\triangle)v\|_{L_z^2(I; H^{\sigma+\delta_1-\f12})}\leq K_\eta\|\nabla v\|_{L_z^2(I; H^{\sigma+\delta_1-\e})},\\
&&\|(T_a+T_A+T_\beta\cdot\nabla)v\|_{L_z^2(I; H^{\sigma+\delta_1-\f12})}\leq K_\eta\|\nabla v\|_{L_z^2(I; L^2)}.
\end{eqnarray*}
By the definition of $A$ and the proof of Lemma \ref{lem:coeffi}, we have
\beno
\|\pa_\xi^\al \pa_zA(z;\cdot,\xi)\|_{C^{-\f12+\e}}\leq K_\eta(1+|\xi|)^{1-|\al|}\quad \textrm{for any} \quad |\xi|\geq \f12,
\eeno
from which and Remark \ref{rem:symb}, it follows that
\beno
\|T_{\pa_z A}v\|_{L_z^2(I; H^{\sigma+\delta_1-\f12})}\leq K_\eta\|\nabla v\|_{L_z^2(I; H^{\sigma+\delta_1-\e})}.
\eeno
The proof is completed.\ef

\begin{lemma}\label{lem:al-V}
It holds that
\beno
&&\|\nabla \al\nabla v\|_{L^1_z(I; H^{-\f12})}+\|\nabla \beta\pa_z v\|_{L^1_z(I; H^{-\f12})}+\|\gamma\pa_z v\|_{L^1_z(I; H^{-\f12})}\\
&&\le K_\eta\|\nabla_{x,z}v\|_{L^2(I\times \R^d)}.
\eeno
\end{lemma}

\no{\bf Proof.}\,We get by the proof of Lemma \ref{lem:PP-Hs} that
\beno
&&\|T_{\gamma}\pa_z v\|_{L^1_z(I; H^{-\f12})}+\|T_{\pa_z v}\gamma\|_{L^1_z(I; H^{-\f12})}\\
&&\le \|\gamma\|_{L^2_z(I;L^\infty)}\|\pa_z v\|_{L^2(I\times \R^d)}\le \|\gamma\|_{\widetilde{L}^2_z(I;C^\e)}\|\pa_z v\|_{L^2(I\times \R^d)},
\eeno
and by the proof of Lemma \ref{lem:PR}, we see that
\begin{align*}
\|R(\gamma, \pa_zv)\|_{L^1_z(I; H^{-\f12})}&\le\int_I\Big(\sum_j2^{-j}\big(\sum_{|k-\ell|
\le 2; k,\ell\ge j-N_0}\|\Delta_k\gamma\|_{L^\infty}\|\Delta_\ell\pa_zv\|_{L^2}\big)^2\Big)^\f12dz\\
&\le C\int_I \sum_j2^{-\f j 2}\sum_{|k-\ell|
\le 2; k,\ell\ge j-N_0}\|\Delta_k\gamma\|_{L^\infty}\|\Delta_\ell\pa_zv\|_{L^2}dz\\
&\le C\|\gamma\|_{\widetilde{L}^2_z(I;C^\e)}\|\pa_z v\|_{L^2(I\times \R^d)},
\end{align*}
which along with Lemma \ref{lem:coeffi} and (\ref{Bony}) give
\beno
\|\gamma\pa_z v\|_{L^1_z(I; H^{-\f12})}\le K_\eta\|\pa_z v\|_{L^2(I\times \R^d)}.
\eeno
The estimates for the other two terms are similar.\ef

Now let us turn to the proof of Proposition \ref{prop:elliptic Hs est}.\vspace{0.1cm}

\no{\bf Proof of Proposition \ref{prop:elliptic Hs est}}. First of all, we consider $\sigma=-\f12$. We have
\beno
\big(\nabla v(z),\nabla v(z)\big)_{H^{-\f12}}
\leq 2\int^z_{-\infty}\big(\pa_{z'}\nabla v(z'),\nabla v(z')\big)_{H^{-\f12}}dz'\le 2\|\nabla_{x,z}v\|^2_{L^2(I\times\R^d)},
\eeno
and by the equation (\ref{eq:ellitic-flat}) and Lemma \ref{lem:al-V}, we get
\begin{align*}
\big(\pa_z {v(z)},\pa_z v(z)\big)_{H^{-\f12}}
&=2\int^z_{-\infty}\big(\pa_{z'}^2 v(z'), \pa_z v(z')\big)_{H^{-\f12}}dz'\\
&=2\int^z_{-\infty}\big(F_0-\al\Delta v+\beta\nabla \pa_zv-\gamma\pa_z v,\pa_z v\big)_{H^{-\f12}}dz'\\
&\leq \Big(\|F_0\|_{Y^{-\f12}(I)}+\|\dv(\al \nabla v+\beta\pa_z v)\|_{L^2_z(I;H^{-1})}\\
&\qquad+\| \nabla\al\nabla v+\nabla\beta \pa_z v+\gamma\pa_z v\|_{L^1_z(I; H^{-\f12})}\Big)\|\pa_zv\|_{X^{-\f12}(I)}\\
&\le K_\eta\big(\|F_0\|_{Y^{-\f12}(I)}+\|\nabla_{x,z}v\|_{L^2(I\times\R^d)}\big)\|\pa_zv\|_{X^{-\f12}(I)}.
\end{align*}
This implies the case of $\sigma=-\f12$.

For general $\sigma$, we use the bootstrap argument. To this end, let us first assume
\begin{eqnarray*}
&&\|\nabla_{x,z}v\|_{X^r(I)}\leq K_\eta\big(\|\na_{x,z}v\|_{L^2(I\times\mathbf{R}^d)}+\|f\|_{H^{r+1}}+\|F_0\|_{Y^{r}(I)}
 +\|\eta\|_{H^{s+\f12}}\|\nabla_{x,z}v\|_{L^\infty(I\times \R^d)}\big).
\end{eqnarray*}
Then we show that the inequality remains true for $r+\delta_1\le s-\f12$ with $\delta_1\le \f12$,
thus the proposition follows since it is true for $r=-\f12$.

Set $w=(\pa_z-T_A)v$. Then $(v,w)$ satisfies the forward and backward parabolic equation respectively:
\beno
&&(\pa_z-T_a)w=F\quad\textrm{ on }I\times\R^d,\quad w|_{z=-\infty}=0,\\
&&(\pa_z-T_A)v=w\quad\textrm{ on }I\times\R^d,\quad v|_{z=0}=f.
\eeno
By Proposition \ref{prop:parabolic estimate} and Lemma \ref{lem:F1}-Lemma \ref{lem:F3}, we infer that
\begin{align}\label{eq:elliptic-Hs-2}
&\|w\|_{X^{r+\delta_1}(I)}\nonumber\\
&\leq K_\eta\big(\|F\|_{Y^{r+\delta_1}(I)}+\|w\|_{L^2_z(I; H^{r+\delta_1})}\big)\nonumber\\
&\le K_\eta \big(\|F_0\|_{Y^{r+\delta_1}(I)}+\|\nabla_{x,z}v\|_{L_z^2(I; H^{r+\delta_1-\e})}
+\|\nabla_{x,z} v\|_{L^\infty(I\times\R^d)}\|\eta\|_{H^{s+\f12}}\big).
\end{align}
Here we use the estimate
\beno
\|w\|_{L^2_z(I; H^{r+\delta_1})}\le K_\eta\|\na_{x,z}v\|_{L^2_z(I; H^{r+\delta_1})}\quad(\textrm{ by Proposition }\ref{prop:symbolic calculus}).
\eeno
Take $\na$ to the equation of $v$ to get
\beno
(\pa_z-T_A)\na v=\na w+T_{\na A}v\quad\textrm{ on }I\times\R^d,\quad \na v|_{z=0}=\na f.
\eeno
By Remark \ref{rem:symb} and Remark \ref{rem:symbol-aA}, we have
\beno
\|T_{\na A}v\|_{L^2_z(I; H^{r+\delta_1-\f12})}\le K_\eta\|\na v\|_{L^2_z(I;H^{r+\delta_1-\e})}.
\eeno
Then by Proposition \ref{prop:parabolic estimate} and (\ref{eq:elliptic-Hs-2}), we get by using $\pa_zv=T_Av+w$ that
\begin{align*}
&\|\na_{x,z}v\|_{X^{r+\delta_1}(I)}\\
&\leq K_\eta\big(\|w\|_{X^{r+\delta_1}(I)}+\|\na v\|_{L^2_z(I;H^{r+\delta_1})}+\|f\|_{H^{r+1+\delta_1}}\big)\\
&\leq K_\eta\big(\|f\|_{H^{r+1+\delta_1}}+\|F_0\|_{Y^{r+\delta_1}(I)}+\|\nabla_{x,z}v\|_{L_z^2(I; H^{r+\delta_1-\e})}
+\|\nabla_{x,z} v\|_{L^\infty(I\times\R^d)}\|\eta\|_{H^{s+\f12}}\big)\\
&\le K_\eta\big(\|\na_{x,z}v\|_{L^2(I\times\R^d)}+\|f\|_{H^{r+\delta_1+1}}+\|F_0\|_{Y^{r+\delta_1}(I)}
 +\|\eta\|_{H^{s+\f12}}\|\nabla_{x,z}v\|_{L^\infty(I\times \R^d)}\big).
\end{align*}
This completes the proof of the proposition. \ef

\subsection{Elliptic estimates in Besov space}

\begin{proposition}\label{prop:elliptic Holder est}
Let $q\in [1,\infty]$ and $v$ be a solution of (\ref{eq:ellitic-flat}) on $I\times \R^d$.
Then for any $r\in [0,\f12], \delta_2>0$, it holds that
\beno
&&\|\na_{x,z} v\|_{\widetilde{L}_z^\infty(I;B^{r}_{\infty,q})}+\|\na_{x,z} v\|_{\widetilde{L}_z^2(I;B^{r+\f12}_{\infty,q})}\\
&&\le  K_\eta\big(\|\na f\|_{B^{r}_{\infty,q}}+\|F_0\|_{\widetilde{Y}^{r}_q(I)}+\|\na_{x,z}v\|_{\widetilde{L}^\infty(I; C^{-\delta_2})}
+\|\nabla_{x,z}v\|_{\widetilde{L}_z^2(I; C^{-\delta_2})}\big),
\eeno
where $\widetilde{Y}^{r}_q(I)\eqdef\widetilde{L}_z^1(I;B^{r}_{\infty,q})+\widetilde{L}_z^2(I;B^{r-\f12}_{\infty,q})$.
\end{proposition}

Let us first present the H\"{o}lder estimates of $F$.

\begin{lemma}\label{lem:F1-2}
For any $r\leq \f12$ and $q\in [1,\infty]$, we have
\beno
\|F_1\|_{\widetilde{L}_z^2(I; B^{r-\f12}_{\infty,q})}
\leq K_\eta\big(\|\nabla_{x,z} v\|_{\widetilde{L}_z^\infty(I; C^{-\f \e2})}+\|\pa_z v\|_{\widetilde{L}_z^\infty(I; B^{r-\f12}_{\infty,q})}\big).
\eeno
\end{lemma}

\no{\bf Proof.}\,By Lemma \ref{lem:PP-Holder} and Lemma \ref{lem:PR}, we have
\begin{align*}
&\|T_\gamma\pa_z v\|_{\widetilde{L}_z^2(I; B^{r-\f12}_{\infty,q})}
\leq C \|\gamma\|_{\widetilde{L}_z^2(I; L^\infty)}\|\pa_z v\|_{\widetilde{L}_z^\infty(I; B^{r-\f12}_{\infty,q})},\\
&\|T_{\pa_z v}\gamma \|_{\widetilde{L}_z^2(I; B^{r-\f12}_{\infty,q})}
       \leq C \|\pa_z v\|_{\widetilde{L}_z^\infty(I; C^{-\f\e 2})}\|\gamma\|_{\widetilde{L}_z^2(I; C^{\e})},\\
&\|R(\gamma,\pa_z v)\|_{\widetilde{L}_z^2(I; B^{r-\f12}_{\infty,q})}
         \leq C \|\pa_z v\|_{\widetilde{L}_z^\infty(I; C^{-\f \e2}})\|\gamma\|_{\widetilde{L}_z^2(I; C^{\e})},
\end{align*}
which along with Lemma \ref{lem:coeffi} and (\ref{Bony}) gives the lemma.\ef

\begin{lemma}\label{lem:F2-2}
For any $r\leq \f12$ and $q\in [1,\infty]$, we have
\beno
\|F_2\|_{\widetilde{L}_z^2(I; B^{r-\f12}_{\infty,q})}\leq K_\eta\|\nabla_{x,z} v\|_{\widetilde{L}_z^\infty(I; C^{-\f\e 2})}.
\eeno
\end{lemma}

\no{\bf Proof.}\,Using (\ref{Bony}), we infer from Lemma \ref{lem:PP-Holder} and Lemma \ref{lem:PR} that
\begin{align*}
\|F_2\|_{\widetilde{L}_z^2(I; B^{r-\f12}_{\infty,q})}
&\leq \|(T_\al-\al)\Delta v\|_{\widetilde{L}_z^2(I; B^{r-\f12}_{\infty,q})}+\|(T_\beta-\beta)\nabla \pa_z v\|_{\widetilde{L}_z^2(I; B^{r-\f12}_{\infty,q})}\\
&\leq C\big(\|\al\|_{\widetilde{L}_z^2(I; C^{r+\f12+\e})}+\|\beta\|_{\widetilde{L}_z^2(I; C^{r+\f12+\e})}\big)\|\nabla_{x,z} v\|_{\widetilde{L}_z^\infty(I; C^{-\f\e2})}\\
&\le K_\eta\|\nabla_{x,z} v\|_{\widetilde{L}_z^\infty(I; C^{-\f\e2})}.
\end{align*}
The proof is finished. \ef

\begin{lemma}\label{lem:F3-2}
For any $r\leq \f12$ and $q\in [1,\infty]$, we have
\beno
\|F_3\|_{\widetilde{L}_z^2(I; B^{r-\f12}_{\infty,q})}\leq K_\eta\|\nabla_{x,z}v\|_{\widetilde{L}_z^\infty(I; C^{-\e })}.
\eeno
\end{lemma}

\no{\bf Proof.}\,From the proof of Lemma \ref{lem:Besov}, we see that
\begin{align*}
&\|(T_aT_A-T_\al \Delta)v\|_{\widetilde{L}_z^2(I; B^{r-\f12}_{\infty,q})}\\
&\leq C\big(\widetilde{M}_0^1(a)\|A\|_{\widetilde{L}^2_z(I; M_{1+\e}^1)}+\|a\|_{\widetilde{L}^2_z(I; M_{1+\e}^1)}\widetilde{M}_0^1(A)\big)
\|\nabla v\|_{\widetilde{L}_z^\infty(I; C^{-\e})}\\
&\le K_\eta\|\nabla v\|_{\widetilde{L}_z^\infty(I; C^{-\e})},
\end{align*}
where we denote
$$\|a\|_{\widetilde{L}^2_z(I; M_{1+\e}^1)}
\eqdef \sup_{|\alpha|\leq d/2+1+\rho}\sup_{|\xi|\geq1/2} \|(1+|\xi|)^{|\alpha|-m}\pa^\alpha_\xi a(z;\cdot,\xi)\|_{\widetilde{L}^2_z(I;C^{1+\e})},
$$
and by Lemma \ref{lem:coeffi}, we have
\beno
\|a\|_{\widetilde{L}^2_z(I; M_{1+\e}^1)}+\|A\|_{\widetilde{L}^2_z(I; M_{1+\e}^1)}\le K_\eta.
\eeno
The estimate for the other parts of $F_3$ is similar.\ef\vspace{0.2cm}

Now let us turn to the proof of Proposition \ref{prop:elliptic Holder est}.\vspace{0.1cm}

\no{\bf Proof of Proposition \ref{prop:elliptic Holder est}.}\,Recall that if we set $w=(\pa_z-T_A)v$, then $(v,w)$ satisfies
\beno
&&(\pa_z-T_a)w=F\quad\textrm{ on }I\times\R^d,\quad w|_{z=-\infty}=0,\\
&&(\pa_z-T_A)v=w\quad\textrm{ on }I\times\R^d,\quad v|_{z=0}=f.
\eeno
By Proposition \ref{prop:parabolic estimate-maximal} and Lemma \ref{lem:F1-2}-Lemma \ref{lem:F3-2}, we deduce that
\begin{align}\label{eq:elliptic-est-w}
&\|w\|_{\widetilde{L}_z^\infty(I;{B^{r}_{\infty,q})}}+\|w\|_{\widetilde{L}_z^2(I;B^{r+\f12}_{\infty,q})}\nonumber\\
&\le K_\eta\big(\|F\|_{\widetilde{Y}^{r}_q(I)}+\|w\|_{\widetilde{L}_z^\infty(I; C^{-\delta_2})}+\|w\|_{\widetilde{L}_z^2(I;C^{-\delta_2})}\big)\nonumber\\
&\le {K_\eta}\big(\|F_0\|_{\widetilde{Y}^{r}_q(I)}+\|\nabla_{x,z} v\|_{\widetilde{L}_z^\infty(I; C^{-\f \e2}\cap B^{r-\f12}_{\infty,q})}+\|\na_{x,z} v\|_{\widetilde{L}_z^2(I;C^{-\delta_2})}\big),
\end{align}
and noting that $(\pa_z-T_A)\na v=\na w+T_{\na A}v$, we get by Proposition \ref{prop:parabolic estimate-maximal} that
\begin{align*}
&\|\na v\|_{\widetilde{L}_z^\infty(I;B^{r}_{\infty,q})}+\|\na v\|_{\widetilde{L}_z^2(I;B^{r+\f12}_{\infty,q})}\\
&\leq K_\eta\big(\|\na f\|_{B^{r}_{\infty,q}}+\|w\|_{\widetilde{L}_z^2(I;B^{r+\f12}_{\infty,q})}+\|\na_{x,z}v\|_{\widetilde{L}_z^2(I;B^{r-\e}_{\infty,q}\cap C^{-\delta_2})}
+\|\na_{x,z}v\|_{\widetilde{L}^\infty_z(I; C^{-\delta_2})}\big)\\
&\le K_\eta\big(\|\na f\|_{B^{r}_{\infty,q}}+\|F_0\|_{\widetilde{Y}^{r}_q(I)}+\|\nabla_{x,z} v\|_{\widetilde{L}_z^\infty(I; C^{-\f \e2}\cap B^{r-\f12}_{\infty,q})}
++\|\na_{x,z}v\|_{\widetilde{L}_z^2(I;B^{r-\e}_{\infty,q}\cap C^{-\delta_2})}\big).
\end{align*}
The estimate for $\pa_z v$ can be deduced by using $\pa_zv=T_Av+w$ and (\ref{eq:elliptic-est-w}). Thus, we obtain
\beno
&&\|\na_{x,z}v\|_{\widetilde{L}_z^\infty(I;B^{r}_{\infty,q})}+\|\na_{x,z}v\|_{\widetilde{L}_z^2(I;B^{r+\f12}_{\infty,q})}\\
&&\le  K_\eta\big(\|\na f\|_{B^{r}_{\infty,q}}+\|F_0\|_{\widetilde{Y}^{r}_q(I)}
+\|\nabla_{x,z} v\|_{\widetilde{L}_z^\infty(I; C^{-\f \e2}\cap B^{r-\f12}_{\infty,q})}+\|\na_{x,z}v\|_{\widetilde{L}_z^2(I;B^{r-\e}_{\infty,q}\cap C^{-\delta_2})}\big),
\eeno
from which and the interpolation, we conclude the proof of the proposition.\ef

\section{Dirichlet-Neumann operator}

\subsection{Definition and basic properties}

We consider the boundary value problem
\ben\label{eq:elliptic-DN}
\left\{
\begin{array}{l}
\Delta_{x,y} \phi=0\qquad\textrm{ in}\quad \cS,\\
\phi|_{y=\eta}=f,
\end{array}
\right.
\een
where $\cS=\big\{(x,y): x\in \R^d, y<\eta(x)\big\}$.
Given $f\in H^\f12(\R^d)$, the existence of the variation solution $\phi$ with $\na_{x,y}\phi\in L^2(\cS)$ can be deduced by using
Riesz theorem, see \cite{ABZ3} for example. Moreover, it holds that
\ben\label{eq:elliptic-H1}
\|\na_{x,y}\phi\|_{L^2(\cS)}\le C(\|\na\eta\|_{L^\infty})\|f\|_{H^\f12}.
\een

\begin{definition}
Given $\eta, f, \phi$ as above, the Dirichlet-Neumann operator $G(\eta)$ is defined by
\beno
G(\eta)f\eqdef\sqrt{1+|\nabla \eta|^2}\pa_n \phi|_{y=\eta}.
\eeno
\end{definition}

We have the following basic properties for $G(\eta)$, see \cite{Lan}.

\begin{proposition}\label{prop:DN-basic properties}
It holds that

\begin{itemize}

\item[1.] the operator $G(\eta)$ is self-adjoint:
\[
(G(\eta)f,g)=(f,G(\eta)g),\quad\forall f,g\in H^\f12(\R^d);
\]

\item[2.] the operator $G(\eta)$ is positive:
\[
(G(\eta)f,f)=\|\na_{x,y}\phi\|_{L^2(\cS)}\ge 0,\quad \forall f\in H^\f12(\R^d);
\]

\item[3.] for any $f,g\in H^\f12(\R^d)$, we have
\begin{eqnarray*}
|(G(\eta)f,g)|\le C(\|\na\eta\|_{L^\infty})\|f\|_{H^\frac12}\|g\|_{H^\frac12};
\end{eqnarray*}

\item[4.] the shape derivative $d_{\eta}G(\eta)$ of $G(\eta)$ is
\beno
d_{\eta}G(\eta)\psi\cdot \pa_t\eta= -G(\eta)(\pa_t\eta B)-\dv(\pa_t\eta V),
\eeno
where $V=\nabla \phi|_{y=\eta}, B=\pa_y \phi|_{y=\eta}.$
\end{itemize}
\end{proposition}

\begin{remark}\label{rem:V-B}
By the definition of Dirichlet-Neumann operator $G(\eta)$, it is easy to see that
\beno
B=\f{\nabla \eta\cdot \nabla \psi+G(\eta)\psi}{1+|\nabla \eta|^2},\quad V=\nabla \psi-B\nabla \eta.
\eeno
\end{remark}

With the notations in Section 4, we denote $v(x,z)=\phi(x,\rho_\delta(x,z))$.
 In terms of $v$, the Dirichlet-Neumann operator $G(\eta)$ can be written as
\ben\label{eq:DN-new def}
G(\eta)f=\Big(\f{1+|\nabla \rho_\delta|^2}{\pa_z \rho_\delta}\pa_z v-\nabla \rho_\delta \cdot \nabla v\Big)\Big|_{z=0}.
\een

\subsection{Tame estimates of the Dirichlet-Neumann operator}
In this subsection, we assume that $\eta\in H^{s+\f12}(\R^d)\cap C^{\f32+\e}(\R^d)$ for $s>1+\f d2$ and $\e>0$.
We denote $I=(-\infty,0)$ and by $K_\eta=K_\eta\big(\|\eta\|_{C^{\f32+\e}}, \|\eta\|_{L^2}\big)$ an increasing function.

Following \cite{ABZ3}, we first paralinearize $G(\eta)$. We set
\beno
\zeta_1=\f{1+|\nabla \rho_\delta|^2}{\pa_z \rho_\delta},\quad \zeta_2=\na\rho_\delta.
\eeno
By Lemma \ref{lem:nonlinear} and (\ref{eq:rho-Holder}), we have
\ben\label{eq:zeta-S-est}
&&\|\zeta_1-1\|_{\widetilde{L}^\infty_z(I;H^{s-\f12})}+\|\zeta_2\|_{\widetilde{L}^\infty_z(I;H^{s-\f12})}\le K_\eta\|\eta\|_{H^{s+\f12}},\\
&&\|\zeta_1\|_{\widetilde{L}^\infty_z(I;C^{\f12+\e})}+\|\zeta_2\|_{\widetilde{L}^\infty_z(I;C^{\f12+\e})}\le K_\eta.\label{eq:zeta-H-est}
\een

Using Bony's decomposition (\ref{Bony}), we decompose $G(\eta)$ as
\beno
G(\eta)f=T_{\zeta_1}\pa_zv+T_{\pa_z v}\zeta_1-R(\zeta_1,\pa_z v)-T_{i\zeta_2\cdot \xi}v-T_{\nabla v}\cdot\zeta_2-R(\zeta_2,\nabla v)\big|_{z=0}.
\eeno
Replacing $\pa_z v$ by $T_Av$, we get
\ben
G(\eta)f=T_{\lambda}f+R(\eta)f,
\een
where $\lambda=\zeta_1A-i\zeta_2\cdot \xi$ with  $A=\f{1}{2}({-i\beta\cdot\xi}+\sqrt{4\al|\xi|^2-(\beta\cdot\xi)^2})$ and
\begin{align*}
R(\eta)f
=&\Big[\big(T_{\zeta_1}T_A-T_{\zeta_1 A}\big)v-T_{\zeta_1}(\pa_z-T_A)v\nonumber\\
  &\quad+\big(T_{\pa_z v}\zeta_1+R(\zeta_1,\pa_z v)-T_{\nabla v}\cdot\zeta_2-R(\nabla v,\zeta_2)\big)\Big]\bigg|_{z=0}\nonumber\\
\equiv & R_1(\eta)f+R_2(\eta)f+R_3(\eta)f.
\end{align*}

\begin{proposition}\label{prop:DN-Hs}
It holds that
\beno
&&\|R(\eta)f\|_{{H^{s-\f12}}}
\leq K_\eta\big(\|f\|_{H^s}+\|\nabla_{x,z}v\|_{L^\infty(I\times\R^d)}\|\eta\|_{H^{s+\f12}}\big),\\
&&\|R(\eta)f\|_{{H^{s-1}}}
\leq  K_\eta\big(\|f\|_{H^{s-\f12}}+\|\nabla_{x,z}v\|_{L^\infty(I\times\R^d)}\|\eta\|_{H^{s+\f12}}\big).
\eeno
\end{proposition}

\no{\bf Proof.}\,Recalling that $A\in \Gamma^1_{\f12+\e}(I\times \R^d)$, we get by Proposition \ref{prop:symbolic calculus}, (\ref{eq:zeta-H-est}), Proposition \ref{prop:elliptic Hs est}
and (\ref{eq:elliptic-H1}) that
\begin{align*}
\|R_1(\eta)f\|_{H^{s-\f12}}\le& K_\eta\|\na v\|_{L^\infty_z(I;H^{s-1})}\\
\le& K_\eta\big(\|\na_{x,z}v\|_{L^2(I\times\R^d)}+\|f\|_{H^s}+\|\nabla_{x,z}v\|_{\widetilde{L}^\infty_z(I;{\color{blue}L^\infty})}\|\eta\|_{H^{s+\f12}}\big)\\
\le& K_\eta\big(\|f\|_{H^s}+\|\nabla_{x,z}v\|_{{L}^\infty(I\times\R^d)}\|\eta\|_{H^{s+\f12}}\big),
\end{align*}
and by Proposition \ref{prop:symbolic calculus}, we have
\begin{align*}
\|R_2(\eta)f\|_{H^{s-\f12}}&\le K_\eta\|(\pa_z-T_A)v\|_{L^\infty_z(I;H^{s-1})}\le K_\eta\|\nabla_{x,z}v\|_{X^{s-1}(I)}\\
&\le K_\eta\big(\|f\|_{H^s}+\|\nabla_{x,z}v\|_{{L}^\infty(I\times\R^d)}\|\eta\|_{H^{s+\f12}}\big).
\end{align*}
For $R_3(\eta)$, we infer from Lemma \ref{lem:PP-Hs}, Lemma \ref{lem:PR} and (\ref{eq:zeta-S-est}) that
\beno
\|R_3(\eta)f\|_{H^{s-\f12}}\le K_\eta\|\nabla_{x,z}v\|_{L^\infty(I\times\R^d)}\|\eta\|_{H^{s+\f12}}.
\eeno
This gives the first inequality. The proof of the second inequality is similar. \ef

\begin{proposition}\label{prop:DN-Holder}
For any $\delta_3>0$, it holds that
\beno
\|R(\eta)f\|_{{C^{\f12}}}\leq K_\eta\big(\|f\|_{C^1}+\|\na_{x,z}v\|_{\widetilde{L}^\infty_z(I;C^{-\delta_3})}+\|\na_{x,z}v\|_{\widetilde{L}^2_z(I;C^{-\delta_3})}\big).
\eeno
\end{proposition}

\no{\bf Proof.}\,By  Lemma \ref{lem:Besov}, Proposition \ref{prop:elliptic Holder est} and (\ref{eq:zeta-H-est}), we get
\beno
\|R_1(\eta)f\|_{C^\f12}\le K_\eta\|\na v\|_{L^\infty_z(I; C^0)}\le K_\eta\big(\|f\|_{C^1}+\|\na_{x,z}v\|_{\widetilde{L}^\infty_z(I;C^{-\delta_3})}
+\|\na_{x,z}v\|_{\widetilde{L}^2_z(I;C^{-\delta_3})}\big).
\eeno
Set $w=(\pa_z-T_A)v$. From the proof of Proposition \ref{prop:elliptic Holder est}(with $r=\f12$), we see that
\beno
\|w\|_{\widetilde{L}^\infty_z(I;C^\f12)}\le K_\eta\|\nabla_{x,z} v\|_{\widetilde{L}_z^\infty(I; C^{0})},
\eeno
hence,
\begin{align*}
\|R_2(\eta)f\|_{C^\f12}&\le K_\eta\|(\pa_z-T_A)v\|_{L^\infty_z(I;C^\f12)}\\
&\le K_\eta\big(\|f\|_{C^1}+\|\na_{x,z}v\|_{\widetilde{L}^\infty_z(I;C^{-\delta_3})}+\|\na_{x,z}v\|_{\widetilde{L}^2_z(I;C^{-\delta_3})}\big).
\end{align*}
It is easy to show by Lemma \ref{lem:PP-Holder} and Lemma \ref{lem:PR} that
 \begin{align*}
\|R_3(\eta)f\|_{C^\f12}
&\le K_\eta\|\na_{x,z} v\|_{L^\infty_z(I; C^0)}\\
&\le K_\eta\big(\|f\|_{C^1}+\|\na_{x,z}v\|_{\widetilde{L}^\infty_z(I;C^{-\delta_3})}+\|\na_{x,z}v\|_{\widetilde{L}^2_z(I;C^{-\delta_3})}\big).
\end{align*}
This finishes the proof.\ef

\begin{remark}
The estimates of $R(\eta)$ may not be optimal, but it is suitable for our application. We refer to
\cite{ABZ3} for more sharper estimates.
\end{remark}

\begin{remark}\label{rem:DN}
By  \eqref{eq:zeta-H-est}  and Remark \ref{rem:symbol-aA}, we know that $\lambda\in \Gamma^1_{\f12+\e}(I\times \R^d)$ with the bound
\beno
\widetilde{M}_{\f12+\e}^1(\lambda)\le K_\eta.
\eeno
\end{remark}

\section{The estimate of the pressure}

Throughout this section,  we denote $\Om_t=\{(x,y): x\in \R^d, y<\eta(t,x)\}, I=(-\infty,0)$, $\widetilde{f}(t,x,z)=f(t,x,\rho_\delta(x,z))$,
and by   $K_\eta=K(\|\eta\|_{C^{\f32+\e}}, \|\eta\|_{L^2})$ a nondecreasing function.

\subsection{The estimates of the velocity potential}
Recall that the velocity potential $\phi$ satisfies
\ben\label{eq:velocity potential}
\Delta_{x,y}\phi=0\quad \hbox{in}\quad \Omega_t,\qquad \phi|_{y=\eta}=\psi.
\een
We infer from Proposition \ref{prop:DN-basic properties} that
\ben\label{eq:phi-H1}
\|\na_{x,y}\phi\|_{L^2(\Om_t)}=\big(G(\eta)\psi,\psi\big)^\f12\equiv E_0(\psi).
\een
By the definition of $(V,B)$, we have
\beno
&&\Delta_{x,y} (\nabla \phi)=0\quad \hbox{in}\quad \Omega_t,
\quad\nabla \phi|_{y=\eta}=V,\\
&&\Delta_{x,y} (\pa_y\phi)=0\quad \hbox{in}\quad \Omega_t,
\quad\pa_y \phi|_{y=\eta}=B.
\eeno
Then it follows from (\ref{eq:elliptic-H1}) and the maximum principle that
\ben
&&\|\nabla_{x,y}^2\phi\|_{L^2(\Om_t))}\leq K_\eta\|(V, B)\|_{H^{\f12}},\label{eq:phi-H2}\\
&&\|\nabla_{x,y}\phi\|_{L^\infty(\Om_t))}\leq \|(V, B)\|_{L^\infty}.\label{eq:phi-G-Linfty}
\een

Further more, we find that
\beno
\left\{\begin{array}{l}
\Delta_{x,y} \pa_{x_i}(\nabla \phi)=0\quad \hbox{in}\quad \Omega_t,\\
\pa_{x_i}\nabla \phi|_{y=\eta}=\pa_{x_i}V-\pa_{x_i}\eta\Big(\na B+\f {\textrm{div} V-\na B\cdot\na \eta} {{1+|\na\eta|^2}}\Big),
\end{array}\right.
\eeno
and
\beno
\left\{\begin{array}{l}
\Delta_{x,y} \pa_{y}(\nabla \phi)=0\quad \hbox{in}\quad \Omega_t,\\
\pa_{y}\nabla \phi|_{y=\eta}=\na B+\na\eta\Big(\f {\textrm{div} V-\na B\cdot\na \eta} {{1+|\na\eta|^2}}\Big),
\end{array}\right.
\eeno
and
\beno
\left\{\begin{array}{l}
\Delta_{x,y}(\pa_{y}^2\phi)=0\quad \hbox{in}\quad \Omega_t,\\
\pa_{y}^2\phi|_{y=\eta}=-\f {\textrm{div} V-\na B\cdot\na \eta} {{1+|\na\eta|^2}}.
\end{array}\right.
\eeno
Using the maximum principle, we deduce that
\ben\label{eq:phi-G2-infty}
\|\na_{x,y}^2\phi\|_{L^\infty(\Om_t)}\le K_\eta\|(\na V, \na B)\|_{L^\infty}.
\een
Then by Proposition \ref{prop:elliptic Hs est} and (\ref{eq:phi-H2}), we infer that
\begin{align}\label{eq:phi-Hs}
\|\nabla_{x,z}(\widetilde{\nabla \phi}, \widetilde{\pa_y\phi})\|_{X^{s-1}(I)}
&\leq K_\eta\big(\|(V, B)\|_{H^{s}}+\|\eta\|_{H^{s+\f12}}\|\nabla^2_{x,y}\phi\|_{L^{\infty}(\Om_t)}\big)\nonumber\\
&\le K_\eta\big(\|(V, B)\|_{H^{s}}+\|\eta\|_{H^{s+\f12}}\|(\na V, \na B)\|_{L^\infty}\big).
\end{align}
And by Proposition \ref{prop:elliptic Holder est} with ${\delta_2>\f d2}$, (\ref{eq:phi-G-Linfty}) and (\ref{eq:phi-H1}), we get
\ben
&&\|\nabla_{x,z}\widetilde{\phi}\|_{\widetilde{L}^\infty_z(I;C^\f12)}
+\|\nabla_{x,z}\widetilde{\phi}\|_{\widetilde{L}^2_z(I;C^1)}\nonumber\\
&&\le K_\eta\big(\|\na\psi\|_{C^\f12}+\|\na_{x,z}\widetilde{\phi}\|_{\widetilde{L}^\infty(I; C^{-\delta_2})}
+\|\nabla_{x,z}\widetilde{\phi}\|_{\widetilde{L}_z^2(I; C^{-\delta_2})}\big)\nonumber\\
&&\le K_\eta\big(\|\na\psi\|_{C^\f12}+\|\nabla_{x,z}\widetilde{\phi}\|_{L^\infty(I\times\R^d)}+\|\nabla_{x,z}\widetilde{\phi}\|_{{L}^2(I\times \R^d)}\big)\nonumber\\
&&\le K_\eta\big(\|(V,B)\|_{W^{1,\infty}}+E_0(\psi)\big).\label{eq:phi-H-est1}
\een
Here we use the fact that $\na \psi=V+B\na \eta$ so that $\|\na \psi\|_{C^\f12}\le K_\eta\|(V,B)\|_{W^{1,\infty}}$.

Using (\ref{Bony}), we get by Lemma \ref{lem:PP-Holder} and Lemma \ref{lem:PR} that
\beno
&&\|\al \triangle \widetilde{\phi}+\beta\cdot\nabla \pa_z \widetilde{\phi}\|_{\widetilde{L}^2_z(I;{ C^0})}
\le C\|\na_{x,z}\widetilde{\phi}\|_{\widetilde{L}^2_z(I;C^1)}\|(\al,\beta)\|_{\widetilde{L}^\infty_z(I;C^\e)},\\
&&\|\gamma \pa_z \widetilde{\phi}\|_{\widetilde{L}^2_z(I;C^0)}\le C\|\pa_z\widetilde{\phi}\|_{\widetilde{L}^\infty_z(I;C^0)}
\|\gamma\|_{\widetilde{L}^2_z(I;C^\e)}.
\eeno
Hence, using the equation
\beno
\pa^2_z \widetilde{\phi}=-\al \triangle \widetilde{\phi}-\beta\cdot\nabla \pa_z \widetilde{\phi}+\gamma \pa_z \widetilde{\phi},
\eeno
we infer from Lemma \ref{lem:coeffi} that
\ben\label{eq:phi-H-est2}
\|\pa_z^2\widetilde{\phi}\|_{\widetilde{L}^2_z(I;C^0)}\le K_\eta\big(\|(V,B)\|_{W^{1,\infty}}+E_0(\psi)\big).
\een
Noticing that
\beno
\pa_z\widetilde{\pa_y\phi}=\pa_z\Big(\f {\pa_z\widetilde{\phi}} {\pa_z\rho_\delta}\Big)
=\f {\pa_z^2\widetilde{\phi}} {\pa_z\rho_\delta}
-\f {\pa_z\widetilde{\phi}\pa_z^2\rho_\delta} {(\pa_z\rho_\delta)^2},
\eeno
thus by (\ref{Bony}), Lemma \ref{lem:PP-Holder}, Lemma \ref{lem:PR}, (\ref{eq:phi-H-est1}), (\ref{eq:phi-H-est2})
and (\ref{eq:rho-H-L2}), we have
\begin{align*}
\|\pa_z\widetilde{\pa_y\phi}\|_{\widetilde{L}_z^2(I; C^{0})}
&\le K_\eta\big(\|\pa_z^2\widetilde{\phi}\|_{\widetilde{L}_z^2(I; C^{0})}
+\|\pa_z\widetilde{\phi}\|_{\widetilde{L}_z^\infty(I; C^{0})}\|\pa_z^2\rho_\delta\|_{\widetilde{L}_z^2(I; C^{\e})}
\big)\\
&\le K_\eta\big(\|(V,B)\|_{W^{1,\infty}}+E_0(\psi)\big).
\end{align*}
Similarly, we can prove
\begin{align*}
\|\na_{x,z}\widetilde{\na_{x,y}\phi}\|_{\widetilde{L}_z^2(I; C^{0})} \le K_\eta\big(\|(V,B)\|_{W^{1,\infty}}+E_0(\psi)\big).
\end{align*}
Then by Proposition \ref{prop:elliptic Holder est} with ${\delta_2>\f d2}$ again, (\ref{eq:phi-G2-infty}) and (\ref{eq:phi-H1}), we get
\ben
&&\|\nabla_{x,z}(\widetilde{\nabla \phi}, \widetilde{\pa_y\phi})\|_{\widetilde{L}^\infty_z(I;C^0)}
+\|\nabla_{x,z}(\widetilde{\nabla \phi}, \widetilde{\pa_y\phi})\|_{\widetilde{L}^2_z(I;C^\f12)}\nonumber\\
&&\le K_\eta\big(\|(V,B)\|_{C^1}+\|\na_{x,z}(\widetilde{\nabla \phi}, \widetilde{\pa_y\phi})\|_{ L^\infty(I\times \R^d)}+\|\nabla_{x,z}(\widetilde{\nabla \phi}, \widetilde{\pa_y\phi})\|_{\widetilde{L}_z^2(I; C^{-\delta_2})}\big)\nonumber\\
&&\le K_\eta\big(\|(V,B)\|_{W^{1,\infty}}+E_0(\psi)\big).\label{eq:phi-H-est}
\een

As an application of (\ref{eq:phi-H-est}), we infer from Proposition \ref{prop:DN-Hs}
and { Proposition \ref{prop:DN-Holder}} that
\begin{lemma}\label{lem:DN-R}
It holds that
\beno
&&\|R(\eta)(V,B)\|_{H^{s-\f12}}
\le K_\eta\big(\|(V,B)\|_{H^s}+(\|(V,B)\|_{W^{1,\infty}}+E_0(\psi))\|\eta\|_{H^{s+\f12}}\big),\\
&&\|R(\eta)(V,B)\|_{H^{s-1}}
\le K_\eta\big(\|(V,B)\|_{H^{s-\f12}}+(\|(V,B)\|_{W^{1,\infty}}+E_0(\psi))\|\eta\|_{H^{s+\f12}}\big),\\
&&\|R(\eta)(V,B)\|_{C^{\f12}}\le K_\eta\big(\|(V,B)\|_{W^{1,\infty}}+E_0(\psi)\big).
\eeno
\end{lemma}

\subsection{The estimates of the pressure}
Recall that the pressure $P$ satisfies
$$
-P=\pa_t\phi+\f{1}{2}|\nabla_{x,y}\phi|^2+gy.
$$
Take $\Delta_{x,y}$ on both sides and use the fact $\Delta_{x,y}\phi=0$ to get
\ben\label{eq:P}
\Delta_{x,y} P=-|\nabla^2_{x,y}\phi|^2=-\na_{x,y}\cdot\big(\na_{x,y}^2\phi\cdot\na_{x,y}\phi\big)\quad \hbox{in}\quad \Omega_t,  \qquad P|_{y=\eta}=0.
\een
By the $L^2$ energy estimate and (\ref{eq:phi-G2-infty}), we get
\begin{align*}
\|\na_{x,y}(P-y)\|_{L^2(\Om_t)}&\le \|\na_{x,y}\phi\na^2_{x,y}\phi\|_{L^2(\Om_t)}\\
&\le K_\eta\|(\widetilde{\na\phi}, \widetilde{\pa_y\phi})\nabla_{x,z}(\widetilde{\nabla \phi}, \widetilde{\pa_y\phi})\|_{L^2(I\times \R^d)}\\
&\le K_\eta\|(\widetilde{\na\phi}, \widetilde{\pa_y\phi})\|_{L^2(I\times \R^d)}
\|\nabla_{x,z}(\widetilde{\nabla \phi}, \widetilde{\pa_y\phi})\|_{L^\infty(I\times \R^d)}\\
&\le K_\eta E_0(\psi)\|(V,B)\|_{W^{1,\infty}},
\end{align*}
which implies that
\ben\label{eq:pressure-H1}
\|\na_{x}\widetilde{P_1}\|_{L^\infty(I;H^{-\f12})}+\|\na_{x,z}\widetilde{P_1}\|_{L^2(I\times \R^d)}\le K_\eta E_0(\psi)\|(V,B)\|_{W^{1,\infty}}.
\een
Here and in what follows, we denote $P_1=P-y$. Indeed, if $y<\eta(x)-2$, following the proof of Proposition \ref{prop:eta-Holde}, we can get
\beno
|\na_{x,y}P_1(x, y)|\le C\big(E_0(\psi)+\|(V,B)\|_{W^{1,\infty}}\big)^2.
\eeno
For $z\in [-2,0]$, we have
\begin{align*}
\pa_z\widetilde{P_1}(x,z)-\pa_z\widetilde{P_1}(x,-2)&=\int_{-2}^z\pa_z^2 \widetilde{P}_1(z')dz'\\
&=\int_{-2}^0\big(F_0-\al\Delta \widetilde{P}_1+\beta\nabla \pa_z\widetilde{P}_1-\gamma\pa_z \widetilde{P}_1\big)dz'
\end{align*}
with $F_0=-\al|\widetilde{\nabla^2_{x,y}\phi}|^2$. Using (\ref{Bony}), we get by Lemma \ref{lem:PP-Holder}, Lemma \ref{lem:PR} and Lemma \ref{lem:coeffi} that
\beno
\|\al\Delta \widetilde{P}_1-\beta\nabla \pa_z\widetilde{P}_1+\gamma\pa_z \widetilde{P}_1\|_{L^1_z(-2,0;H^{-1})}\le K_\eta\|\na_{x,z}\widetilde{P_1}\|_{L^2}
\eeno
So, we conclude that
\ben\label{eq:pressure-infty}
\|\pa_z \widetilde{P_1}\|_{L^\infty_z(I;L^\infty+H^{-1})}\le K_\eta\big(E_0(\psi)+\|(V,B)\|_{W^{1,\infty}}\big)^2.
\een

\begin{lemma}\label{lem:P-source1}
It holds that
\beno
&&\|\al \widetilde{|\nabla^2_{x,y}\phi|^2}\|_{\widetilde{L}_z^2(I; B^{0}_{\infty,1})}
\le K_\eta\big(E_0(\psi)+\|(V,B)\|_{W^{1,\infty}}\big)^2,\\
&&\|\al \widetilde{|\nabla^2_{x,y}\phi|^2}\|_{\widetilde{L}_z^2(I; H^{s-1})}
\le K_\eta\big(\|(V,B)\|_{W^{1,\infty}}\|(V,B)\|_{H^s}+\|(V,B)\|_{W^{1,\infty}}^2\|\eta\|_{H^{s+\f12}}\big).
\eeno
\end{lemma}

\no{\bf Proof.}\,Using (\ref{Bony}), we get by Lemma \ref{lem:PP-Holder} and
Lemma \ref{lem:PR} along with Lemma \ref{lem:coeffi} and (\ref{eq:rho-Holder}) that
\begin{align*}
\|\al \widetilde{|\nabla^2_{x,y}\phi|^2}\|_{\widetilde{L}_z^2(I; B^{0}_{\infty,1})}
&\le K_\eta\|\widetilde{|\nabla^2_{x,y}\phi|^2}\|_{\widetilde{L}_z^2(I; B^{0}_{\infty,1})}\\
&\le K_\eta\|\widetilde{\nabla^2_{x,y}\phi}\|_{\widetilde{L}_z^\infty(I; C^{0})}
\|\widetilde{\nabla^2_{x,y}\phi}\|_{\widetilde{L}_z^2(I; C^{\f12})}\\
&\le K_\eta\|\na_{x,z}(\widetilde{\nabla_{x,y}\phi})\|_{\widetilde{L}_z^\infty(I; C^0)}
\|\na_{x,z}(\widetilde{\nabla_{x,y}\phi})\|_{\widetilde{L}_z^2(I; C^{\f12})},
\end{align*}
where we use the chain rule so that
\beno
\widetilde{\pa_{x_i}\na_{x,y}\phi}=\pa_{x_i}\widetilde{\na_{x,y}\phi}-\pa_{x_i}\rho_\delta\f {\pa_z\widetilde{\na_{x,y}\phi}} {\pa_z\rho_\delta},
\quad \widetilde{\pa_y\na_{x,y}\phi}=\f {\pa_z\widetilde{\na_{x,y}\phi}} {\pa_z\rho_\delta}.
\eeno
Then the first inequality of the lemma follows from (\ref{eq:phi-G2-infty}) and (\ref{eq:phi-H-est}).

Using (\ref{Bony}), we get by Lemma \ref{lem:PP-Hs}-
Lemma \ref{lem:PR} along with Lemma \ref{lem:coeffi} and (\ref{eq:rho-Hs})-(\ref{eq:rho-Holder}) that
\begin{align*}
&\|\al \widetilde{|\nabla^2_{x,y}\phi|^2}\|_{\widetilde{L}_z^2(I; H^{s-1})}\\
&\le K_\eta\big(\|\na^2_{x,y}\phi\|^2_{L^\infty(\Om_t)} \|\eta\|_{H^{s+\f12}}
+\|\na^2_{x,y}\phi\|_{L^\infty(\Om_t)} \|\widetilde{\nabla^2_{x,y}\phi}\|_{L^2_z(I;H^{s-1})}\big)\\
&\le K_\eta\big(\|\na^2_{x,y}\phi\|^2_{L^\infty(\Om_t)} \|\eta\|_{H^{s+\f12}}
+\|\na^2_{x,y}\phi\|_{L^\infty(\Om_t)} \|\na_{x,z}\widetilde{\nabla_{x,y}\phi}\|_{L^2_z(I;H^{s-1})}\big),
\end{align*}
from which, (\ref{eq:phi-Hs}) and (\ref{eq:phi-G2-infty}), we deduce the second inequality.\ef
\vspace{0.1cm}

Now we infer from Proposition \ref{prop:elliptic Holder est} with $\delta_2>\f {d}2+1$, Lemma \ref{lem:P-source1}, (\ref{eq:pressure-H1})
and (\ref{eq:pressure-infty}) that
\ben\label{eq:pressure-Holder}
&&\|\nabla_{x,z}\widetilde{P_1}\|_{\widetilde{L}_z^\infty(I; B^{\f12}_{\infty,1})}+\|\nabla_{x,z} \widetilde{P_1}\|_{\widetilde{L}_z^2(I; B^{1}_{\infty,1})}\nonumber\\
&&\leq K_\eta\big(\|\al \widetilde{|\nabla^2_{x,y}\phi|^2}\|_{\widetilde{L}_z^2(I; B^{0}_{\infty,1})}
   +\|\nabla_{x,z}\widetilde{P_1}\|_{\widetilde{L}^\infty_z(I; C^{-\delta_2})}
   +\|\nabla_{x,z}\widetilde{P_1}\|_{\widetilde{L}^2_z(I;C^{-\delta_2})}\big)\nonumber\\
&&\leq K_\eta\big(\|\al \widetilde{|\nabla^2_{x,y}\phi|^2}\|_{\widetilde{L}_z^2(I; B^{0}_{\infty,1})}
   +\|\nabla_{x,z}\widetilde{P_1}\|_{L^\infty_z(I; L^\infty+H^{-1})} +\|\nabla_{x,z}\widetilde{P_1}\|_{L^2(I\times \R^d)}\big)\nonumber\\
&&\le K_\eta\big(E_0(\psi)+\|(V,B)\|_{W^{1,\infty}}\big)^2.
\een
On the other hand, using the equation
\ben\label{eq:P-z}
\pa^2_z \widetilde{P_1}=-\al \triangle \widetilde{P_1}-\beta\cdot\nabla \pa_z \widetilde{P_1}+\gamma \pa_z \widetilde{P_1}-\al \widetilde{|\nabla^2_{x,y}\phi|^2},
\een
we infer from Lemma \ref{lem:PP-Holder}, Lemma \ref{lem:PR}, Lemma \ref{lem:coeffi} and Lemma \ref{lem:P-source1} that
\begin{align}\label{eq:pressure-H-2}
\|\pa_z^2\widetilde{P_1}\|_{\widetilde{L}_z^2(I; B^{0}_{\infty,1})}
&\le C\|(\al,\beta)\|_{\widetilde{L}_z^\infty(I; C^{\e})}
\|\na_{x,z}\widetilde{P_1}\|_{\widetilde{L}_z^2(I; B^{1}_{\infty,1})}\nonumber\\
&\qquad+C\|\gamma\|_{\widetilde{L}_z^2(I; C^{\e})}\|\pa_z\widetilde{P_1}\|_{\widetilde{L}_z^\infty(I; B^{\f12}_{\infty,1})}
+\|\al \widetilde{|\nabla^2_{x,y}\phi|^2}\|_{\widetilde{L}_z^2(I; B^{0}_{\infty,1})}\nonumber\\
&\le K_\eta\big(E_0(\psi)+\|(V,B)\|_{W^{1,\infty}}\big)^2.
\end{align}

While, it follows from Proposition \ref{prop:elliptic Hs est}, Lemma \ref{lem:P-source1},
(\ref{eq:pressure-H1}) and (\ref{eq:P-z}) that
\begin{eqnarray}
&&\|\nabla_{x,z} \widetilde{P_1}\|_{X^{s-\f12}(I)}+\|\pa_z^2\widetilde{P_1}\|_{L^2_z(I;H^{s-1})}\nonumber\\
&&\leq K_\eta\big(\|\na_{x,z}\widetilde{P_1}\|_{L^2(I\times \R^d)}+ \|\al \widetilde{|\nabla^2_{x,y}\phi|^2}\|_{\widetilde{L}_z^2(I; H^{s-1})}
   +\|\eta\|_{H^{s+\f12}}\|\nabla_{x,z} \widetilde{P_1}\|_{L^\infty(I\times \R^d))} \big)\nonumber\\
&&\le K_\eta\big(1+E_0(\psi)+\|(V,B)\|_{W^{1,\infty}}\big)^2\big(\|\eta\|_{H^{s+\f12}}+\|(V,B)\|_{H^s}\big).\label{eq:pressure-Hs}
\end{eqnarray}
Here we use the fact $\|\nabla_{x,z} \widetilde{P_1}\|_{L^\infty(I\times \R^d)} \le C \|\nabla_{x,z}\widetilde{P_1}\|_{\widetilde{L}_z^\infty(I; B^{\f12}_{\infty,1})}$
and
\beno
\|\pa_z^2\widetilde{P_1}\|_{L^2_z(I;H^{s-1})}\le K_\eta\big(\|\na_{x,z}\widetilde{P_1}\|_{L^\infty(I\times \R^d)}\|\eta\|_{H^{s+\f12}}+\|\nabla_{x,z} \widetilde{P_1}\|_{X^{s-\f12}(I)}+
\|\al \widetilde{|\nabla^2_{x,y}\phi|^2}\|_{\widetilde{L}_z^2(I; H^{s-1})}\big),
\eeno
which follows from Lemma \ref{lem:PP-Hs}, Lemma \ref{lem:PR} and Lemma \ref{lem:coeffi}.
\vspace{0.1cm}

To estimate $(\pa_t+V\cdot\nabla)a$, we derive the equation of  $\dot{P}\eqdef (\pa_t+\nabla_{x,y} \phi\cdot\na_{x,y})P$.

\begin{lemma}\label{lem:P}
Assume that $(\phi, \eta, P)$ is a smooth solution of the water-wave system (\ref{eq:vp})-(\ref{eq:euler-ber}). Then we have
\ben\label{eq:P-D}
\left\{
\begin{array}{l}
\Delta_{x,y} \dot{P}=4\nabla_{x,y}^2\phi\cdot\nabla_{x,y}^2 P+2\sum_{i,j,k}(\pa_i\pa_j\phi)(\pa_i\pa_k\phi)(\pa_j\pa_k\phi)\equiv F\quad \hbox{in}\quad \Omega_t,\\
\dot{P}|_{y=\eta}=0.
\end{array}
\right.
\een
\end{lemma}

\no{\bf Proof.}\,By (\ref{eq:vp}) and (\ref{eq:P}), we get
\begin{align*}
\Delta_{x,y}(\na_{x,y}\phi\cdot\na_{x,y}P)&=2\na^2_{x,y}\phi\cdot\na_{x,y}^2P+\na_{x,y}\phi\cdot\Delta_{x,y}\na_{x,y}P\\
&=2\na^2_{x,y}\phi\cdot\na_{x,y}^2P-2\pa_k\phi(\pa_i\pa_j\phi)(\pa_i\pa_j\pa_k\phi).
\end{align*}
Hence,
\beno
\Delta_{x,y} \dot{P}=2\nabla_{x,y}^2\phi\cdot\nabla_{x,y}^2 P-2(\pa_i\pa_j\phi)\big(\pa_t(\pa_i\pa_j\phi)+\pa_k\phi\pa_k(\pa_i\pa_j\phi)\big).
\eeno
Taking $\pa_i\pa_j$ on both sides of (\ref{eq:euler-ber}), we get
\beno
\pa_t(\pa_i\pa_j\phi)+\pa_k\phi\pa_k(\pa_i\pa_j\phi)=-\pa_i\pa_jP-(\pa_i\pa_k\phi)(\pa_j\pa_k\phi).
\eeno
This gives the first equation.

Due to $P(t,x,\eta)=0$, we infer that
\beno
P_t+\eta_t\pa_yP|_{y=\eta}=0, \quad \na P+\na\eta\pa_y P|_{y=\eta}=0,
\eeno
which implies that
\beno
P_t+\na_{x,y}\phi\cdot\na_{x,y} P|_{y=\eta}=-\pa_y P(\eta_t+\na\phi\cdot\na \eta-\pa_y\phi)|_{y=\eta},
\eeno
from which and (\ref{eq:vp-b}), it follows that $\dot P|_{y=\eta}=0$.\ef

\begin{remark}\label{rem:div form}
Using $\Delta_{x,y}\phi=0$, the second term on the right hand side of (\ref{eq:P-D}) can be written as the
divergence form:
\beno
&&\sum_{i,j,k}(\pa_i\pa_j\phi)(\pa_i\pa_k\phi)(\pa_j\pa_k\phi)\\
&&=\sum_{i,j,k}\pa_i\big((\pa_j\phi)(\pa_i\pa_k\phi)(\pa_j\pa_k\phi)\big)-\sum_{i,j,k}(\pa_j\phi)(\pa_i\pa_k\phi)(\pa_i\pa_j\pa_k\phi)\\
&&=\sum_{i,j,k}\pa_i\big((\pa_j\phi)(\pa_i\pa_k\phi)(\pa_j\pa_k\phi)\big)-\f12\sum_{i,j,k}(\pa_j\phi)\pa_j\big((\pa_i\pa_k\phi)(\pa_i\pa_k\phi)\big)\\
&&=\sum_{i,j,k}\pa_i\big((\pa_j\phi)(\pa_i\pa_k\phi)(\pa_j\pa_k\phi)\big)-\f12\na_{x,y}\cdot\big(|\na_{x,y}^2\phi|^2\na\phi\big).
\eeno
\end{remark}

Now we infer from (\ref{eq:P-D}) and Remark \ref{rem:div form} that
\beno
\|\na_{x,y}\dot P\|_{L^2(\Om_t)}\le 4\|\na^2_{x,y}\phi\|_{L^\infty(\Om_t)}\|\na_{x,y}P_1\|_{L^2(\Om_t)}+
2\|\na^2_{x,y}\phi\|_{L^\infty(\Om_t)}^2\|\na_{x,y}\phi\|_{L^2(\Om_t)}.
\eeno
Then following the proof of (\ref{eq:pressure-infty}), we get by (\ref{eq:elliptic-H1}), (\ref{eq:phi-H1}), (\ref{eq:phi-G2-infty}), (\ref{eq:pressure-H1})
and (\ref{eq:pressure-Holder}) that
\ben\label{eq:P-D-H1}
\|\na_{x,z}\widetilde{\dot P}\|_{L^\infty(I;L^\infty+H^{-1})}+\|\na_{x,z}\widetilde{\dot P}\|_{L^2(I\times \R^d)}\le K_\eta\big(E_0(\psi)+\|(V,B)\|_{W^{1,\infty}}\big)^3.
\een

\begin{lemma}\label{lem:P-non}
It hods that
\beno
&&\|\al\widetilde{F}\|_{\widetilde{L}_z^1(I; B^0_{\infty,1})}\le K_\eta\big(\|(V,B)\|_{W^{1,\infty}}+E_0(\psi)\big)^3.
\eeno
\end{lemma}

\no{\bf Proof.}\,Using (\ref{Bony}), we get by Lemma \ref{lem:PP-Holder} and Lemma \ref{lem:PR} that
\begin{align*}
\|\al\widetilde{F}\|_{\widetilde{L}_z^1(I; B^0_{\infty,1})}&\le K_\eta\big(\|\widetilde{\na_{x,y}^2\phi}\cdot\widetilde{\na_{x,y}^2P}\|_{\widetilde{L}_z^1(I; B^0_{\infty,1})}
+\|\widetilde{(\na_{x,y}^2\phi)^2}\|_{\widetilde{L}_z^1(I; C^\e)}\|\na^2_{x,y}\phi\|_{L^\infty(\Om_t)}\big)\\
&\le K_\eta\big(\|\widetilde{\na_{x,y}^2\phi}\|_{\widetilde{L}_z^2(I; C^\e)}\|\widetilde{\na_{x,y}^2P}\|_{\widetilde{L}_z^2(I; B^0_{\infty,1})}
+\|\widetilde{\na_{x,y}^2\phi}\|_{\widetilde{L}_z^2(I; C^\e)}^2\|\na^2_{x,y}\phi\|_{L^\infty(\Om_t)}\big).
\end{align*}
By the chain rule, we have that for $i=1,\cdots,d$,
\beno
&&\pa_i\widetilde{P}=\widetilde{\pa_iP}+\widetilde{\pa_zP}\cdot\pa_i\rho_\delta,\quad\pa_z\widetilde{P}=\widetilde{\pa_zP}\cdot\pa_z\rho_\delta,\quad (\pa_i=\pa_{x_i})\\
&&\pa_z^2\widetilde{P}=\widetilde{\pa_y^2P}\cdot(\pa_z\rho_\delta)^2+\widetilde{\pa_y P}\cdot\pa_z^2\rho_\delta,\\
&&\pa^2_{i,j}\widetilde{P}=\widetilde{\pa^2_{i,j}P}+\widetilde{\pa^2_{i,z}P}\cdot \pa_j\rho_\delta+\widetilde{\pa^2_{j,z}P}\cdot \pa_i\rho_\delta +\pa_{i,j}^2\rho_\delta\cdot\widetilde{\pa_z P}+\pa_i\rho_\delta\pa_j\rho_\delta\cdot \widetilde{\pa^2_{zz}P}, \\
&&\pa^2_{i,z}\widetilde{P}=\widetilde{\pa^2_{i,z}P}\cdot\pa_z\rho_\delta+\pa^2_{i,z}\rho_\delta\cdot\widetilde{\pa_iP}
+\pa_z\rho_\delta\pa_i\rho_\delta\cdot\widetilde{\pa^2_{i,z}P},
\eeno
which along with Lemma \ref{lem:PP-Holder}, Lemma \ref{lem:PR}, (\ref{eq:rho-Holder}) and (\ref{eq:rho-H-L2}) implies that
\beno
\|\widetilde{\na_{x,y}^2P}\|_{\widetilde{L}_z^2(I; B^0_{\infty,1})}
\le K_\eta\big(\|\nabla_{x,z}\widetilde{P_1}\|_{\widetilde{L}_z^\infty(I; B^{\f12}_{\infty,1})}+\|\pa_z^2\widetilde{P_1}\|_{\widetilde{L}_z^2(I; B^{0}_{\infty,1})}\big).
\eeno
Then the lemma follows from (\ref{eq:phi-G-Linfty}), (\ref{eq:phi-H-est}), (\ref{eq:pressure-Holder})
and (\ref{eq:pressure-H-2}).\ef\vspace{0.1cm}

Then we infer from Proposition \ref{prop:elliptic Holder est} with $\delta_2>\f {d}2+1$, Lemma \ref{lem:P-non}
and (\ref{eq:P-D-H1}) that
\begin{align}
\|\na_{x,z}\widetilde{\dot P}\|_{\widetilde{L}_z^\infty(I; B^0_{\infty,1})}&\le K_\eta\big(\|\al\widetilde{F}\|_{\widetilde{L}_z^1(I; B^0_{\infty,1})}+
\|\nabla_{x,z}\widetilde{\dot P}\|_{L^\infty_z(I; L^\infty+H^{-1})} +\|\nabla_{x,z}\widetilde{\dot P}\|_{L^2(I\times \R^d)}\big)\nonumber\\
&\le K_\eta\big(\|(V,B)\|_{W^{1,\infty}}+E_0(\psi)\big)^3.\label{eq:P-D-Holder}
\end{align}

\section{New formulation and symmetrization}

Recall the water-wave system
\beq \label{eq:euler-zak-7}
\left\{
\begin{array}{ll}
\pa_t \eta-G(\eta)\psi=0, \\
\pa_t \psi+g\eta+\frac{1}{2}|\na
\psi|^2-\frac{(G(\eta)\psi+\na\eta\cdot\na
\psi)^2}{2(1+|\na\eta|^2)}=0.
\end{array}\right.
\eeq
Following the framework of \cite{ABZ3}, we introduce the new unknowns
\ben
\zeta=\nabla \eta ,\quad B=\pa_y\phi|_{y=\eta},\quad V=\nabla\phi|_{y=\eta},\quad a=-\pa_y P|_{y=\eta}.
\een
Recall that the pressure $P$ satisfies
\ben\label{eq:pressure}
-P=\pa_t\phi+\f{1}{2}|\nabla_{x,y}\phi|^2+gy,
\een
where $\phi$ is the solution of the elliptic equation
\beno
\Delta_{x,y}\phi=0\quad \hbox{in}\quad \Omega_t,\qquad \phi|_{y=\eta}=\psi.
\eeno

Then the system (\ref{eq:euler-zak-7}) can be reformulated as(see \cite{ABZ3}):

\begin{lemma}\label{prop:new formulation}
The new unknowns $(V, B, \zeta)$ satisfy
\ben
&&(\pa_t+V\cdot\nabla)B=a-g,\label{B}\label{eq:B}\\
&&(\pa_t+V\cdot\nabla)V+a\zeta=0,\label{V}\label{eq:V}\\
&&(\pa_t+V\cdot\nabla)\zeta=G(\eta)V+\zeta G(\eta)B\label{eq:zeta}.
\een
\end{lemma}

\no{\bf Proof.} For the reader's convenience, we present a proof. By the chain rule,
for any function $f=f(t,x,y)$, we have
\begin{align}
(\pa_t+V\cdot\nabla)(f|_{y=\eta})
&=(\pa_t+V\cdot\nabla)(f(t,x,\eta))\nonumber\\
&=\big[(\pa_t f+\nabla\phi\cdot\nabla f)+\pa_y f(\pa_t\eta+V\cdot\nabla\eta)\big]\big|_{y=\eta}\nonumber\\
&=(\pa_tf+\nabla_{x,y}\phi\cdot\nabla_{x,y}f)|_{y=\eta}.\label{eq:chain}
\end{align}
Here in the last equality we use the fact that
\ben
\pa_t\eta+V\cdot\nabla\eta=B.\label{eq:eta}
\een

Taking $\na_{x,y}$ to (\ref{eq:pressure}), we deduce the equalities (\ref{eq:B})-(\ref{eq:V}) from (\ref{eq:chain}) and $P(t,x,\eta)=0$.
Taking $\pa_{x_i}$ to (\ref{eq:eta}), we get
\ben
(\pa_t+V\cdot\nabla)\pa_{x_i}\eta=\pa_{x_i}B-\sum_{j}\pa_{x_i}V_j\pa_{x_j}\eta.\label{eq:eta-derivative}
\een
By the definitions of $(V,B)$ and $G(\eta)$, we find
\begin{align*}
\pa_{x_i}B-\sum_{j}\pa_{x_i}V_j\pa_{x_j}\eta
=&\pa_y\pa_{x_i}\phi-\sum_{j}\pa_{x_j}\eta\pa_{x_j}\pa_{x_i}\phi\big|_{y=\eta}
  +\pa_{x_i}\eta\big(\pa_y(\pa_y\phi)-\nabla\eta\cdot\nabla\pa_y\phi\big)\big|_{y=\eta}\\
=&G(\eta)V_i+\pa_{x_i}\eta G(\eta)B,
\end{align*}
which along with (\ref{eq:eta-derivative}) gives (\ref{eq:zeta}).\ef

Now we introduce the so called good unknown $(U_s, \zeta_s)$ defined by
\beno
U_s\eqdef \langle D\rangle^s V+ T_\zeta \langle D\rangle^s B,\quad \zeta_s\eqdef \langle D\rangle^s\zeta.
\eeno

\begin{lemma}
The unknown $(U_s, \zeta_s)$ satisfies
\ben\label{eq:good unknown}
\left\{
\begin{array}{l}
(\pa_t+T_V\cdot \nabla)U_s+T_a\zeta_s=f_1,\\
(\pa_t+T_V\cdot \nabla)\zeta_s=T_\lambda U_s+f_2,
\end{array}
\right.
\een
where $(f_1, f_2)$ is given by
\begin{align*}
f_1=&\Ds h_1-[\Ds, T_V\cdot \nabla]V-[\Ds,T_a]\zeta-[\Ds,T_\zeta](\pa_t+T_V\cdot \nabla)B\nonumber\\
    &-T_\zeta[\Ds,T_V\cdot \nabla]B-[T_\zeta,\pa_t+T_V\cdot \nabla]\Ds B,\\
f_2=&\Ds h_2-[\Ds, T_V\cdot \nabla]\zeta+[T_\lambda,\Ds]U\nonumber,
\end{align*}
with $U=V+T_\zeta B$ and $(h_1,h_2)$ given by
\beno
&& h_1=(T_V-V)\cdot\nabla V-R(a,\zeta)+T_\zeta(T_V-V)\cdot\nabla B,\\
&&h_2=(T_V-V)\cdot\na \zeta+[T_\zeta, T_\lambda]B+(\zeta-T_\zeta)T_\lambda B+R(\eta)V+\zeta R(\eta)B.
\eeno
\end{lemma}

\no{\bf Proof.}\,Applying  Bony's decomposition (\ref{Bony})  to (\ref{eq:B})-(\ref{eq:zeta}),  we get
\beno
&&(\pa_t+T_V\cdot\na)V+T_a\zeta+T_\zeta(\pa_t+T_V\cdot\na)B=h_1,\\
&&(\pa_t+T_V\cdot\na)\zeta=T_\lambda U+h_2.
\eeno
Then the system (\ref{eq:good unknown}) follows by applying  $\Ds$ to the above equations. \ef\vspace{0.1cm}

We denote
\beno
\gamma=\sqrt{a\lambda},\quad q=\sqrt{\f{a}{\lambda}},\quad \theta_s=T_q \zeta_s.
\eeno
Taking $T_q$ on the both sides of the second equation of (\ref{eq:good unknown}),  we obtain the following symmetrized system:
\ben\label{eq:symmetrized system}
\left\{
\begin{array}{l}
(\pa_t+T_V\cdot \nabla)U_s+T_\gamma\theta_s=F_1,\\
(\pa_t+T_V\cdot \nabla)\theta_s-T_\gamma U_s=F_2,
\end{array}
\right.
\een
with $(F_1,F_2)$ given by
\beno
&&F_1=f_1+(T_\gamma T_q-T_a)\zeta_s,\\
&&F_2=T_q f_2+(T_qT_\lambda-T_\gamma)U_s-[T_q,\pa_t+T_V\cdot \nabla]\zeta_s.
\eeno

\section{Energy estimates}

Assume that $(U_s, \theta_s)$ is a solution of  (\ref{eq:symmetrized system}) on $[0,T]$, and $a(t,x)$ satisfies
\ben
\inf_{(t,x)\in [0,T]\times \R^d}a(t,x)\ge c_0.
\een
We denote by $K_\eta^1=K_\eta^1\Big(\sup_{t\in [0,T]}\big(\|\eta(t)\|_{C^{\f32+\e}}+\|\eta(t)\|_{L^2}\big), c_0^{-1}\Big)$ by an increasing function,
which may be different from line to line. By the definition of $(\gamma,q)$, it is easy to show that
\ben\label{eq:gammaq-est}
 M^{\f12}_{0}(\gamma)+ M^{-\f12}_0(q)\le K_\eta^1\|a\|^{\f12}_{C^\f12},\quad M^{\f12}_{\f12}(\gamma)+ M^{-\f12}_{\f12}(q)\le K_\eta^1\|a\|_{C^\f12}.
\een
And by Lemma \ref{lem:PP-Hs}, we have
\beno
\|(U_s,\theta_s)\|_{L^2}\le  K_\eta^1\big(\|(V,B)\|_{H^s}+\|\eta\|_{H^{s+\f12}}\big).
\eeno

The goal of this section is to prove that

\begin{proposition}
\ben\label{eq:energy-Utheta}
\f{d}{dt}\|(U_s,\theta_s)\|_{L^2}\le K_\eta^1\big(G_1(t)\|(V,B)\|_{H^s}+G_2(t)\|\eta\|_{H^{s+\f12}}+\|a-g\|_{H^{s-\f12}}\big),
\een
where $G_i(i=1,2)$ is defined by
\beno
&&G_1(t)=1+\|a\|_{C^{\f12}}+\|(V, B)\|_{B^1_{\infty, 1}},\\
&&G_2(t)=1+\|a\|^{\f32}_{C^{\f12}}+\|\pa_t a+V\cdot\nabla a\|_{L^\infty}+\|a\|_{C^{\f12}}\big(\|(V, B)\|_{B^1_{\infty, 1}}
+E_0(\psi)\big).
\eeno
\end{proposition}

\no{\bf Proof.}\,We multiply $(U_s, \theta_s)$ by both sides of (\ref{eq:symmetrized system}) and integrate on $\mathbf{R}^d$ to obtain
\ben\label{eq:energy-1}
\f12\f{d}{dt}\|(U_s,\theta_s)\|^2_{L^2}=I_1+I_2+I_3,
\een
with $I_i$ given by
\beno
&&I_1=-\langle T_V\cdot \nabla U_s,U_s\rangle-\langle T_V\cdot\nabla \theta_s,\theta_s\rangle,\nonumber\\
&&I_2=-\langle T_\gamma\theta_s,U_s\rangle+\langle T_\gamma U_s,\theta_s\rangle,\nonumber\\
&&I_3=\langle F_1,U_s\rangle+\langle F_s,\theta_s\rangle.
\eeno
By Proposition \ref{prop:symbolic calculus}, we know that
\beno
&&\|(T_V\cdot \nabla)^*+T_V\cdot \nabla\|_{L^2\rightarrow L^2}\leq C\|V\|_{W^{1,\infty}},\nonumber\\
&&\|(T_\gamma-(T_\gamma)^*\|_{L^2\rightarrow L^2}\leq C M^{\f12}_{\f12}(\gamma)\leq K_\eta^1\|a\|_{C^{\f12}},
\eeno
from which and (\ref{eq:energy-1}), we infer that
\ben\label{eq:energy-2}
\f{d}{dt}\|(U_s,\theta_s)\|_{L^2}
\leq K_\eta^1\big(\|(V,B)\|_{W^{1,\infty}}+\|a\|_{C^{\f12}}\big)\|(U_s,\theta_s)\|_{L^2}+\|(F_1,F_2)\|_{L^2}.
\een

It remains to estimate $\|(F_1,F_2)\|_{L^2}$.
By Proposition \ref{prop:symbolic calculus} and (\ref{eq:gammaq-est}), we get
\begin{align*}
&\|(T_\gamma T_q-T_a)\zeta_s\|_{L^2}
\leq K_\eta^1\|a\|_{C^\f12}^\f32\|\zeta_s\|_{H^{-\f12}}\le K_\eta^1\|a\|_{C^\f12}^\f32\|\eta\|_{H^{s+\f12}},\\
&\|(T_qT_\lambda-T_\gamma)U_s\|_{L^2}
\leq K_\eta^1\|a\|_{C^\f12}\|U_s\|_{L^2}\le K_\eta^1\|a\|_{C^\f12}\|(V,B)\|_{H^{s}}.
\end{align*}
By Proposition \ref{prop:symbolic calculus} and Proposition \ref{prop:commu-Ds}, we get
\beno
&&\|[\Ds,T_V\cdot \nabla]V\|_{L^2}\leq C\|V\|_{W^{1,\infty}}\|V\|_{H^{s}},\\
&&\|T_\zeta[\Ds, T_V\cdot \nabla]B\|_{L^2}\le K_\eta^1\|V\|_{W^{1,\infty}}\|B\|_{H^s},\\
&&\|T_q[\Ds,T_V\cdot \nabla]\zeta\|_{L^2}
\leq K(\eta, a)\|a\|_{C^\f12}\|V\|_{W^{1,\infty}}\|\eta\|_{H^{s+\f12}},
\eeno
and by Remark \ref{rem:DN},
\beno
&&\|[\Ds,T_a]\zeta\|_{L^2}\leq C\|a\|_{C^{\f12}}\|\na\zeta\|_{H^{s-\f12}}\le C\|a\|_{C^{\f12}}\|\na\zeta\|_{H^{s-\f12}},\\
&&\|T_q[\Ds,T_\lambda]U\|_{L^2}\leq K_\eta^1\|a\|_{C^{\f12}}\|U\|_{H^{s}}\le  K_\eta^1\|a\|_{C^{\f12}}\|(V,B)\|_{H^{s}}.
\eeno

Using the equation $\pa_t B+V\cdot\na B=a-g$, we get by Proposition \ref{prop:commu-Ds} that
\begin{align*}
\|[\Ds, T_\zeta](\pa_tB+T_V\cdot\nabla B)\|_{L^2}
&\leq K_\eta^1\|\pa_tB+T_V\cdot\nabla B\|_{H^{s-\f12}}\nonumber\\
&\leq K_\eta^1\big(\|a-g\|_{H^{s-\f12}}+\|B\|_{W^{1,\infty}}\|V\|_{H^s}\big).
\end{align*}

By Proposition \ref{prop:commutator}, we get
\beno
&&\|[T_\zeta,\pa_t+T_V\cdot \nabla]\Ds B\|_{L^2}\\
&&\leq K_\eta^1\big(\|V\|_{B^{1}_{\infty,1}}+\|\pa_t \zeta +V\cdot\nabla \zeta\|_{L^\infty}\big)\|B\|_{H^s}\nonumber\\
&&\leq K_\eta^1\big(\|V\|_{B^{1}_{\infty,1}}+\|(V,B)\|_{W^{1,\infty}}\big)\|B\|_{H^s}.
\eeno
Here we use the fact that
\beno
(\pa_t+V\cdot\nabla)\pa_i \eta=\pa_i B-\sum_{j=1}^d\pa_i V_j \pa_j\eta.
\eeno
Similarly, we have
\beno
&&\|[T_q,\pa_t+T_V\cdot \nabla]\zeta_s\|_{L^2}\\
&&\leq  C\big(\mathcal{M}^{-\f12}_{0}(q)\|V\|_{B^1_{\infty,1}}+\mathcal{M}^{-\f12}_{0}(\pa_tq+V\cdot\nabla q)\big)\|\zeta_s\|_{H^{-1/2}}\nonumber\\
&&\leq K_\eta^1\big(\|a\|_{C^\f12}\|V\|_{B^1_{\infty,1}}+
\|\pa_t a+V\cdot\nabla a\|_{L^\infty}+\|a\|_{C^\f12}\|\pa_t \nabla \eta+V\cdot\nabla \nabla \eta\|_{L^\infty}\big)\|\eta\|_{H^{s+\f12}}\\
&&\le K_\eta^1\big(\|a\|_{C^\f12}\|V\|_{B^1_{\infty,1}}+
\|\pa_t a+V\cdot\nabla a\|_{L^\infty}+\|a\|_{C^\f12}\|(V,B)\|_{W^{1,\infty}}\big)\|\eta\|_{H^{s+\f12}}.
\eeno

Using (\ref{Bony}), we infer from Lemma \ref{lem:PP-Hs} and Lemma \ref{lem:PR} that
\beno
\|\Ds h_1\|_{L^2}
\leq K_\eta^1\big(\|(V,B)\|_{W^{1,\infty}}\|(V,B)\|_{H^s}+\|a\|_{C^{\f12}}\|\eta\|_{H^{s+\f12}}\big).
\eeno
Next we present the estimate of $h_2$. First of all, we have by Lemma \ref{lem:PP-Hs} and Lemma \ref{lem:PR} that
\beno
&&\|(T_V-V)\cdot\na \zeta\|_{H^{s-\f12}}\le K_\eta^1\|V\|_{H^s},\\
&&\|(\zeta-T_\zeta)T_\lambda B\|_{H^{s-\f12}}\le C\|T_\lambda B\|_{L^\infty}\|\eta\|_{H^{s+\f12}}.
\eeno
It follows from Proposition \ref{prop:symbolic calculus} that
\beno
\|[T_\zeta,T_\lambda]B\|_{H^{s-\f12}}\le K_\eta^1\|B\|_{H^s}.
\eeno
And by Lemma \ref{lem:DN-R}, we get
\begin{align*}
\|R(\eta)V\|_{H^{s-\f12}}&\le K_\eta^1\big(\|V\|_{H^s}+(\|(V,B)\|_{W^{1,\infty}}+E_0(\psi))\|\eta\|_{H^{s+\f12}}\big),\\
\|\zeta R(\eta)B\|_{H^{s-\f12}}&\leq K_\eta^1\big(\|R(\eta)B\|_{H^{s-\f12}}+\|R(\eta)B\|_{L^\infty}\|\zeta\|_{H^{s-\f12}})\nonumber\\
&\leq K_\eta^1\big(\|B\|_{H^s}+(\|(V, B)\|_{W^{1,\infty}}+E_0(\psi))\|\eta\|_{H^{s+\f12}}\big).
\end{align*}
Hence, we deduce that
\beno
&&\|\Ds T_q h_2\|_{L^2}\\
&&\le K_\eta^1\|a\|_{C^\f12}\big(\|(V,B)\|_{H^s}+(\|(V, B)\|_{W^{1,\infty}}+E_0(\psi)+\|T_\lambda B\|_{L^\infty})\|\eta\|_{H^{s+\f12}}\big).
\eeno
Noting that
\beno
\|T_\lambda B\|_{L^\infty}\le \|T_\lambda B\|_{B^0_{\infty,1}}\le K_\eta^1\|B\|_{B^1_{\infty,1}},
\eeno
then by summing up the above estimates, we conclude that
\beno
\|(F_1, F_2)\|_{L^2}\leq K^1_{\eta}\big(G_1(t)\|(V, B)\|_{H^s}+G_2(t)\|\eta\|_{H^{s+\f12}}+\|a-g\|_{H^{s-\f12}}\big),
\eeno
from which and (\ref{eq:energy-2}), we deduce (\ref{eq:energy-Utheta}).\ef
\vspace{0.1cm}

Next we recover the estimate of $(V,B,\eta)$ from that of $(U_s,\theta_s)$.

\begin{lemma}\label{lem:recover}
It holds that
\begin{align*}
&\|\eta\|_{H^{s+\f12}}\le K_\eta^1\big(\|(\eta,\theta_s)\|_{L^2}+\|a\|_{C^\f12}^\f32\|\zeta_s\|_{H^{-1}}\big),\\
&\|(V,B)\|_{H^s}\le K_\eta^1\big(\|U_s\|_{L^2}+\|(V,B)\|_{L^2}+(\|(V,B)\|_{W^{1,\infty}}+E_0(\psi))\|\eta\|_{H^{s+\f12}}\big).
\end{align*}
\end{lemma}

\no{\bf Proof.}\,First of all, we have
\beno
\zeta_s=(1-T_{1/q}T_q)\zeta_s+T_{1/q}T_{q}\zeta_s=(1-T_{1/q}T_q)\zeta_s+T_{1/q}\theta_s,
\eeno
which along with Proposition \ref{prop:symbolic calculus} and (\ref{eq:gammaq-est}) implies
\begin{align*}
\|\zeta_s\|_{H^{-\f12}}&\le K_\eta^1\big(\|a\|_{C^\f12}^\f32\|\zeta_s\|_{H^{-1}}+\|\theta_s\|_{L^2}\big).
\end{align*}
Hence, we get
\ben\label{eq:eta-Hs}
\|\eta\|_{H^{s+\f12}}
\leq \|\eta\|_{L^2}+\|\zeta_s\|_{H^{-\f12}}\leq K_\eta^1\big(\|\eta\|_{L^2}+\|a\|_{C^\f12}^\f32\|\zeta_s\|_{H^{-1}}+\|\theta_s\|_{L^2}\big).
\een

Recall that
\beno
U=V+T_\zeta B,\quad \nabla B=G(\eta)V.
\eeno
Thus, we get
\begin{align*}
\nabla U =& \nabla V+T_\zeta \nabla B+T_{\nabla \zeta} B\\
=&\nabla V+T_\zeta G(\eta)V+T_{\nabla \zeta} B\\
=&\nabla V+T_\zeta(T_\lambda V+R(\eta)V)+T_{\nabla \zeta} B.
\end{align*}
Let $T_\mathfrak{p}V=T_{i\xi+\zeta\lambda}V$, then
\begin{align*}
T_\mathfrak{p}V=& \nabla U-(T_\zeta T_\lambda -T_{\zeta\lambda})V-T_\zeta(R(\eta)V)+(T_{i\xi}V-\nabla V)-T_{\nabla \zeta}B\\
\equiv& \nabla U +R'(\eta)V-T_{\nabla \zeta}B,
\end{align*}
which implies
\ben\label{eq:V-U}
 V=T_{1/\mathfrak{p}}\big( \nabla U +R'(\eta)V-T_{\nabla \zeta}B\big)+(1-T_{1/\mathfrak{p}}T_{\mathfrak{p}})V.
\een
Then by Proposition \ref{prop:symbolic calculus} and Lemma \ref{lem:DN-R}, we get
\begin{align*}
\|V\|_{H^s}\le& K_\eta^1\big(\|U\|_{H^s}+\|R(\eta)V\|_{H^{s-1}}+\|B\|_{H^{s-\f12}}+\|V\|_{H^{s-\f12}}\big)\\
\le& K_\eta^1\big(\|U\|_{H^s}+\|(V,B)\|_{L^2}\big)+\f12\|(V,B)\|_{H^s}\\
&+K_\eta^1(\|(V,B)\|_{W^{1,\infty}}+E_0(\psi))\|\eta\|_{H^{s+\f12}},
\end{align*}
on the other hand, we have
\begin{align*}
 \|B\|_{H^s}&\le \|B\|_{L^2}+\|G(\eta)V\|_{H^{s-1}}\\
 &\le \|B\|_{L^2}+K_\eta^1\big(\|V\|_{H^s}+(\|(V,B)\|_{W^{1,\infty}}+E_0(\psi))\|\eta\|_{H^{s+\f12}}\big),
\end{align*}
from which, it follows that
\ben\label{eq:VB-Hs}
\|(V,B)\|_{H^s}\le K_\eta^1\big(\|U_s\|_{L^2}+\|(V,B)\|_{L^2}+(\|(V,B)\|_{W^{1,\infty}}+E_0(\psi))\|\eta\|_{H^{s+\f12}}\big).
\een
Then the lemma follows from (\ref{eq:eta-Hs}) and (\ref{eq:VB-Hs}). \ef

\begin{lemma}\label{lem:VB-infinity}
It holds that
\beno
\|(V,B)\|_{B^{1}_{\infty,1}}\le  K_\eta^1\big(1+\|(V,B)\|_{W^{1,\infty}}+E_0(\psi)\big)\ln(e+\|U_s\|_{L^2}+\|B\|_{H^s}).
\eeno
\end{lemma}

\no{\bf Proof.}\,It follows from (\ref{eq:V-U}), Lemma \ref{lem:Besov} and Lemma \ref{lem:DN-R} that
\begin{align*}
\|V\|_{B^1_{\infty,1}}&\le K_\eta^1\big(\|U\|_{B^1_{\infty,1}}+\|R(\eta)V\|_{B^0_{\infty,1}}+\|(V,B)\|_{B^\f12_{\infty,1}}\big)\\
&\le  K_\eta^1\big(\|U\|_{B^1_{\infty,1}}+\|(V,B)\|_{W^{1,\infty}}+E_0(\psi)\big),
\end{align*}
hence by $\na B=G(\eta)V$, we get
\begin{align*}
\|\na B\|_{B^0_{\infty,1}}&\le \|G(\eta)V\|_{B^0_{\infty,1}}\le \|T_\lambda V\|_{B^{0}_{\infty,1}}+\|R(\eta)V\|_{C^\f12}\\
&\le K_\eta^1\big(\|V\|_{B^1_{\infty,1}}+\|(V,B)\|_{W^{1,\infty}}+E_0(\psi)\big)\\
&\le  K_\eta^1\big(\|U\|_{B^1_{\infty,1}}+\|(V,B)\|_{W^{1,\infty}}+E_0(\psi)\big).
\end{align*}
This proves that
\beno
\|(V,B)\|_{B^{1}_{\infty,1}}\le  K_\eta^1\big(\|U\|_{B^1_{\infty,1}}+\|(V,B)\|_{W^{1,\infty}}+E_0(\psi)\big).
\eeno
Give any $N\in\textbf{N}$, we have
\begin{align*}
\|U\|_{B^1_{1,\infty}}&\le \sum_{j\le N}2^j\|\Delta_j U\|_{L^\infty}+\sum_{j>N}2^j\|\Delta_jU\|_{L^\infty}\\
&\le CN\|(V,B)\|_{W^{1,\infty}}+\sum_{j>N}2^{(1-s)j}\|\Delta_j\Ds U\|_{L^\infty}\\
&\le  CN\|(V,B)\|_{W^{1,\infty}}+\sum_{j>N}2^{(1-s)j}\|\Delta_j\Ds U\|_{L^\infty}\\
&\le  CN\|(V,B)\|_{W^{1,\infty}}+2^{-N(s-1-\f d2)}\|\Ds U\|_{L^2},
\end{align*}
taking $N$ such that $2^{-N(s-1-\f d2)}\|\Ds U\|_{L^2}\sim 1$, we get
\begin{align*}
\|U\|_{B^1_{1,\infty}}
&\le C \big(1+\|(V,B)\|_{W^{1,\infty}}\big)\ln(e+\| \Ds U\|_{L^2})\\
&\le  K_\eta^1 \big(1+\|(V,B)\|_{W^{1,\infty}}\big)\ln(e+\| U_s\|_{L^2}+\|B\|_{H^s}).
\end{align*}
The proof is finished.\ef

\section{Proof of Theorem \ref{thm:blow-up}}

\subsection{The basic energy law}
We introduce the total energy functional $H(\eta,\psi)$ as
\ben
H(\eta,\psi)\eqdef \int_{\mathbf{R}^d}\big(g|\eta|^2+\psi G(\eta)\psi\big)dx.
\een

\begin{proposition}\label{prop:energy law}
Assume that $(\eta,\psi)$ is a smooth solution of (\ref{eq:euler-zak}) with the initial data $(\eta_0,\psi_0)$
on $[0,T]$. Then it holds that
\beno
H(\eta(t),\psi(t))=H(\eta_0,\psi_0)\qquad \textrm{for any} \quad t\in [0,T].
\eeno
\end{proposition}

\no{\bf Proof.}\,Multiplying $g\eta$ and $G(\eta)\psi$ on the both side of the first and second equation of water wave respectively, and integrating on $\mathbf{R}^d$, then we add the resulting equations to obtain
\beno
(\pa_t\eta, g\eta)+(\pa_t\psi, G(\eta)\psi)=-\f12\Big(|\na \psi|^2-\frac{(G(\eta)\psi+\na\eta\cdot\na \psi)^2}{(1+|\na\eta|^2)},G(\eta)\psi\Big).
\eeno
First of all, we have
\beno
(\pa_t\psi, G(\eta)\psi)=\pa_t(\psi, G(\eta)\psi)-(d_{\eta}G(\eta)\psi\cdot \pa_t\eta, \psi)-(\psi, G(\eta)\pa_t\psi).
\eeno
Then by Proposition \ref{prop:DN-basic properties} and (\ref{eq:euler-zak}), we get
\begin{align}
2(\pa_t\psi, G(\eta)\psi)
&=\pa_t(\psi, G(\eta)\psi)-(d_{\eta}G(\eta)\psi\cdot \pa_t\eta, \psi)\nonumber\\
&=\pa_t(\psi, G(\eta)\psi)+(G(\eta)(\pa_t\eta B)+\dv(\pa_t\eta V), \psi)\nonumber\\
&=\pa_t(\psi, G(\eta)\psi)+(\pa_t\eta B,G(\eta) \psi)-(\pa_t\eta V, \nabla \psi)\nonumber\\
&=\pa_t(\psi, G(\eta)\psi)+(G(\eta) \psi  B,G(\eta) \psi)-(G(\eta) \psi V, \nabla \psi)\nonumber\\
&=\pa_t(\psi, G(\eta)\psi)+(G(\eta) \psi, G(\eta) \psi  B)-(G(\eta) \psi , V\cdot \nabla \psi)\nonumber\\
&=\pa_t(\psi, G(\eta)\psi)+(G(\eta) \psi, G(\eta) \psi  B- V\cdot \nabla \psi).\label{eq:energy law-1}
\end{align}
It follows from Remark \ref{rem:V-B} that
\begin{align*}
|\na \psi|^2-\frac{(G(\eta)\psi+\na\eta\cdot\na \psi)^2}{(1+|\na\eta|^2)}
&=(V+B\nabla\eta)\cdot \nabla \psi-B(G(\eta)\psi+\na\eta\cdot\na \psi)\\
&=V\cdot \nabla \psi-BG(\eta)\psi,
\end{align*}
which along with (\ref{eq:energy law-1}) implies that
\beno
\f {d} {dt}\int_{\mathbf{R}^d}(g|\eta|^2+\psi G(\eta)\psi)dx=0.
\eeno
This implies the proposition.\ef

\subsection{H\"{o}lder estimate of the free surface from the mean curvature}
Recall that the equation of the mean curvature
\begin{eqnarray}\label{eq:mean curvature}
\nabla\cdot \Big(\f{\nabla \eta}{\sqrt{1+|\nabla \eta|^2}}\Big)=\kappa,
\end{eqnarray}
where $\kappa$ is the mean curvature of the free surface $y=\eta(t,x)$.

\begin{proposition}\label{prop:eta-Holde}
Assume that $\eta\in W^{1,\infty}$ and $\kappa\in L^2\cap L^p$ for some $p>d$. Then
$\eta\in C^{2-\f dp}$ with the bound
\beno
\|\eta\|_{C^{2-\f d p}}\le C\big(\|\na\eta\|_{L^\infty}, \|\kappa\|_{L^2\cap L^p}\big).
\eeno
\end{proposition}

\no{\bf Proof.}\,Taking $\pa_\ell=\pa_{x_\ell}$ on both sides of (\ref{eq:mean curvature}), we get
\begin{eqnarray*}
\pa_\ell \kappa=\pa_i\Big((1+|\na\eta|^2)^{-\f12}\pa_i\pa_\ell\eta-(1+|\na\eta|^2)^{-\f32}\pa_j\eta\pa_i\eta\pa_i\pa_\ell\eta\Big).
\eeno
We set $\eta_\ell=\pa_\ell\eta$ and $a_{ij}=(1+|\na\eta|^2)^{-\f32}\big((1+|\na\eta|^2)\delta_{ij}-\pa_i\eta\pa_j\eta\big)$.
Then we find that
\ben
\pa_j\big(a_{ij}\pa_i\eta_\ell\big)=\pa_\ell \kappa.
\een
It is easy to verify that the matrix $\big(a_{ij}\big)$ is uniformly elliptic with the elliptic constants depending on $\|\na\eta\|_{L^\infty}$.
Using the De Giorgi method, it can be proved that $\eta_\ell \in C^\epsilon$ for some $\epsilon>0$ and
\beno
\|\eta_\ell\|_{C^{\epsilon}}\le C(\|\na \eta\|_{L^\infty}, \|\kappa\|_{L^p}).
\eeno
This means that $\eta\in C^{1+\epsilon}$, hence  $a_{ij}\in C^{\epsilon}$.

Next we prove H\"{o}lder regularity of $\eta_\ell$ by freezing the leading coefficients method.
For any ball $B_r(x_0)\subset \R^d$ with radius $r$ and center $x_0$, let $w$ be a unique solution of the Dirichlet problem
\beno
\int_{B_r(x_0)}a_{ij}(x_0)\pa_iw\pa_j\varphi dx=0\quad \textrm{for any }\varphi\in H^1_0(B_r(x_0))
\eeno
with $w-\eta_\ell \in H^1_0(B_r(x_0))$. Then $v=\eta_\ell-w$ satisfies
\beno
\int_{B_r(x_0)}a_{ij}(x_0)\pa_iv\pa_j\varphi dx=\int_{B_r(x_0)}\big(-\kappa\pa_\ell \varphi+(a_{ij}(x_0)-a_{ij}(x))\pa_i\eta_\ell \pa_j\varphi\big) dx
\eeno
for any $\varphi\in H^1_0(B_r(x_0))$. Take $\varphi=v$ to get
\beno
\int_{B_r(x_0)}|\na v|^2dx\le C\Big(r^{2\epsilon}\int_{B_r(x_0)}|\na \eta_\ell|^2dx+\int_{B_r(x_0)}|k|^2dx\Big),
\eeno
which along with Lemma \ref{lem:compare-Harmonic} gives for any $0<\rho\le r$
\beno
\int_{B_\rho(x_0)}|\na \eta_\ell|^2dx \le C\Big(\big(r^{2\epsilon}+\big(\f {\rho} r\big)^{d}\big)\int_{B_r(x_0)}|\na \eta_\ell|^2dx+r^{d(1-\f 2p)}\|\kappa\|_{L^p}^2\Big),
\eeno
from which  and a standard iteration,  we infer that there exists $R_0>0$ such that  for any $0<\rho<r\le R_0$,
\beno
\int_{B_\rho(x_0)}|\na \eta_{\ell}|^2dx\le C\Big( \big(\f {\rho} r\big)^{d(1-\f2p)}\int_{B_r(x_0)}|\na \eta_\ell|^2dx+\rho^{d(1-\f2p)}\|\kappa\|_{L^p}^2\Big).
\eeno
In particular, taking $r=R_0$ yields that for any $\rho<R_0$,
\beno
\int_{B_\rho(x_0)}|\na \eta_\ell|^2dx\le C\rho^{d(1-\f2p)}\big(\|\na \eta_\ell\|_{L^2}^2+\|\kappa\|_{L^p}^2\big)
\le  C\rho^{d(1-\f2p)}\big(\|\kappa\|_{L^2}^2+\|\kappa\|_{L^p}^2\big).
\eeno
This implies that $\eta_\ell\in C^{1-\f d p}$ and
\beno
\|\eta_\ell\|_{C^{1-\f dp}}\le C\big(\|\na \eta\|_{L^\infty}, \|\kappa\|_{L^2\cap L^p}\big).
\eeno
Hence, $\eta\in C^{2-\f d p}$ and the proposition follows.\ef

\begin{lemma}\label{lem:compare-Harmonic}
Let  $w$ be as in the proof of Proposition \ref{prop:eta-Holde}. Then for any $u\in H^1(B_r(x_0))$ and  $0<\rho\le r$,  it hods that
\beno
\int_{B_\rho(x_0)}|\na u|^2dx\le C\Big(\big(\f {\rho} r\big)^d\int_{B_r(x_0)}|\na u|^2dx+\int_{B_r(x_0)}|\na(u-w)|^2dx\Big),
\eeno
where $C$ is a  constant depending only on the elliptic constants of $(a_{ij}(x_0))$.
\end{lemma}

\no{\bf Proof.}\,Set $v=u-w$,  we have for any $0<\rho\le r$,
\begin{align*}
\int_{B_\rho(x_0)}|\na u|^2dx&\le 2\int_{B_\rho(x_0)}|\na w|^2dx+2\int_{B_\rho(x_0)}|\na v|^2dx\\
&\le C\big(\f {\rho} r\big)^d\int_{B_r(x_0)}|\na w|^2dx+2\int_{B_r(x_0)}|\na v|^2dx\\
&\le C\big(\f {\rho} r\big)^d\int_{B_r(x_0)}|\na u|^2dx+C\int_{B_r(x_0)}|\na v|^2dx.
\end{align*}
Here we used the property of Harmonic function for $w$,  since it satisfies an elliptic equation with constant coefficients. \ef

\subsection{Proof of Theorem \ref{thm:blow-up}}

Recall the assumption of the theorem:
\beno
&M(T)\eqdef\sup_{t\in [0,T]}\|\kappa(t)\|_{L^p\cap L^2}+\int_0^T\|(\na V,\na B)(t)\|_{L^\infty}^6dt<+\infty,\\
&\displaystyle\inf_{(t,x,y)\in [0,T]\times \Sigma_t}-\frac {\pa P} {\pa \textbf{n}}(t,x,y)\ge c_0.
\eeno
In order to prove Theorem \ref{thm:blow-up}, it suffices to show that
\beno
\sup_{t\in [0,T]}E_s(t)\le C\big(E_s(0), M(T),T,\textrm{ TS}(a)^{-1}\big),
\eeno
where $C(\cdots)$ is an increasing function, and
\beno
E_s(t)\eqdef\|(\eta,\psi)(t)\|_{H^{s+\f12}}+\|(V,B)(t)\|_{H^s},\quad \textrm{TS}(a)\eqdef \inf_{(t,x)\in [0,T)\times \Om_t}a(t,x).
\eeno
In what  follows, we denote by $K_\eta^1=K_\eta^1\big( \sup_{t\in [0,T]}\big(\|\eta(t)\|_{C^{\f32+\e}}+\|\eta(t)\|_{L^2}\big), \textrm{TS}(a)^{-1}\big)$
an increasing function.

Thanks to the equation of $\eta$, we find that
\beno
(\pa_t+V\cdot\nabla )\pa_i\eta=\pa_iB-\sum_{j=1}^d\pa_iV\pa_j\eta.
\eeno
which implies that
\beno
\|\nabla\eta\|_{L^\infty((0,T)\times \R^d)}\leq C\big(E_s(0), M(T),T\big).
\eeno
Hence by Proposition \ref{prop:energy law} and Proposition \ref{prop:eta-Holde}, we obtain
\beno
\sup_{t\in [0,T]}\big(\|\eta(t)\|_{C^{\f32+\e}}+\|\eta(t)\|_{L^2}\big)\le C\big(E_s(0), M(T),T\big),
\eeno
with $\e=\f12-\f d p>0$. Hence, $K_\eta\le C\big(E_s(0), M(T),T\big)$. Note that
\beno
a(t,x)=-\pa_yP|_{y=\eta}=-\f {1} {\sqrt{1+|\na\eta|^2}}\f {\pa P} {\pa \textbf{n}}|_{y=\eta}.
\eeno
Hence, $\textrm{TS}(a)\ge c_1$ for some $c_1>0$.

By the definition of $a$, we infer from (\ref{eq:pressure-Holder}) and (\ref{eq:pressure-Hs}) that
\ben
&&\|a\|_{C^\f12}\le K_\eta^1\big(\|(V,B)\|_{W^{1,\infty}}+E_0(\psi)\big)^2,\label{eq:a-Holder}\\
&&\|a-g\|_{H^{s-\f12}}\le K_\eta^1\big(1+\|(V,B)\|_{W^{1,\infty}}+E_0(\psi)\big)^2\big(\|\eta\|_{H^{s+\f12}}+\|(V,B)\|_{H^s}\big).\label{eq:a-Hs}
\een
By  (\ref{eq:chain}), we find  that
\beno
\pa_t a+V\cdot\na a=\pa_y\dot P-\pa_y\na_{x,y}\phi\cdot\na_{x,y}P|_{y=\eta},
\eeno
which  along with (\ref{eq:P-D-Holder}), (\ref{eq:pressure-Holder}) and \eqref{eq:phi-G2-infty}  implies
\ben
\|\pa_t a+V\cdot\na a\|_{L^\infty}\le  K_\eta^1\big(\|(V,B)\|_{W^{1,\infty}}+E_0(\psi)\big)^3.\label{eq:a-D-est}
\een

Recall that $(V,B,\zeta)$ satisfies
\beno
&&(\pa_t+V\cdot\nabla)B=a-g,\\
&&(\pa_t+V\cdot\nabla)V+a\zeta=0,\\
&&(\pa_t+V\cdot\nabla)\zeta=G(\eta)V+\zeta G(\eta)B.
\eeno
Making $L^2$ energy estimate for $(V,B)$, we get
\beno
\f{d}{dt}\|(V,B)\|_{L^2}\le \|(\na V,\na B)\|_{L^\infty}\|(V,B)\|_{L^2}+\|a-g\|_{L^2}+\|a\|_{L^\infty}\|\zeta\|_{L^2}.
\eeno
While, making $H^{s-1}$ energy estimates for $\zeta$, it is easy to obtain
\begin{align}
\f d {dt}\|\zeta\|_{H^{s-1}}&\le C\|\na V\|_{L^\infty}\|\zeta\|_{H^{s-1}}
+\|\zeta\|_{L^\infty}\|V\|_{H^s}+\|G(\eta)V\|_{H^{s-1}}+\|\zeta G(\eta)B\|_{H^{s-1}}\nonumber\\
&\le K_\eta^1\big(\|\na V\|_{L^\infty}\|\zeta\|_{H^{s-1}}+\|(V,B)\|_{H^s}+(\|(V,B)\|_{W^{1,\infty}}+E_0(\psi))\|\eta\|_{H^{s+\f12}}\big).\nonumber
\end{align}
Then by (\ref{eq:energy-Utheta}), (\ref{eq:a-Holder})-(\ref{eq:a-D-est}), Lemma \ref{lem:recover} and Lemma \ref{lem:VB-infinity},
we obtain
\beno
&&\f{d}{dt}\big(\|(U_s,\theta_s)\|_{L^2}+\|(V,B)\|_{L^2}+\|\zeta_s\|_{H^{-1}}\big)\\
&&\le K_\eta^1G(t)\big(\|(U_s,\theta_s)\|_{L^2}+\|\zeta_s\|_{H^{-1}}+\|(V,B,\eta)\|_{L^2}\big)
\ln\big(e+\|(U_s,\theta_s)\|_{L^2}+\|(V,B,\eta)\|_{L^2}\big).
\eeno
with $G(t)=\big(1+E_0(\psi)+\|(V,B)\|_{W^{1,\infty}}\big)^6$.
Note that
\beno
&&\|(V,B)\|_{L^\infty}\le C\big(\|(V,B)\|_{H^{-\f12}}+\|(\na V,\na B)\|_{L^\infty}\big)\le
K_\eta\big(E_0(\psi)+\|(\na V,\na B)\|_{L^\infty}\big),\\
&&E_0(\psi)+\|\eta\|_{L^2}\le C\big(\|\psi_0\|_{H^\f12}+\|\eta_0\|_{L^2}\big).
\eeno
Then we apply Gronwall's inequality to obtain
\beno
\|(U_s,\theta_s)\|_{L^2}+\|(V,B)\|_{L^2}\le C\big(E_s(0), M(T),T,\textrm{ TS}(a)^{-1}\big).
\eeno
Noting that for any $\epsilon>0$,
\beno
\|(V,B)\|_{W^{1,\infty}}\le C_\epsilon\|(V,B)\|_{L^2}+\epsilon\|(V,B)\|_{H^s}.
\eeno
Then by Lemma \ref{lem:recover} again, we deduce  that
\beno
\|\eta\|_{H^{s+\f12}}+\|(V,B)\|_{H^s}\le C\big(E_s(0), M(T),T,\textrm{ TS}(a)^{-1}\big).
\eeno
This completes the proof of the theorem.\ef

\section*{Acknowledgements}

{\small The second author thanks Professor Sijue Wu for helpful discussions.
Zhifei Zhang is partly supported by NSF of China under Grant 10990013 and 11071007,
 Program for New Century Excellent Talents in University and Fok Ying Tung Education Foundation.}

\end{document}